\documentclass{amsart}
\usepackage{amssymb}

\newtheorem{thm}{Theorem}[section]

\newtheorem{lem}[thm]{Lemma}
\newtheorem{rema}[thm]{Remark}

\theoremstyle{remark}

\numberwithin{equation}{section}

\newcommand{\h}[0]{\mathfrak{h}}

\newcommand{\Lnu}[0]{\hat{L}_\nu}

\def \Z{\mathbb{Z}}

\def \Q{\mathbb{Q}}
\def \C{\mathbb{C}}

\begin{document}

\title[An equivalence of twisted modules]{An equivalence of two
constructions of permutation-twisted modules for lattice vertex
operator algebras}

\author{Katrina Barron}
\address{Department of Mathematics, University of Notre Dame,
Notre Dame, IN 46556}
\email{kbarron@nd.edu}

\author{Yi-Zhi Huang}
\address{Department of Mathematics, Rutgers University,
New Brunswick, NJ 08854}
\email{yzhuang@math.rutgers.edu}

\author{James Lepowsky}
\address{Department of Mathematics, Rutgers University,
New Brunswick, NJ 08854}
\email{lepowsky@math.rutgers.edu}

\begin{abstract}
The problem of constructing twisted modules for a vertex operator
algebra and an automorphism has been solved in particular in two
contexts.  One of these two constructions is that initiated by the
third author in the case of a lattice vertex operator algebra and an
automorphism arising {}from an arbitrary lattice isometry.  This
construction, {}from a physical point of view, is related to the
space-time geometry associated with the lattice in the sense of string
theory.  The other construction is due to the first author, jointly
with C. Dong and G.  Mason, in the case of a multi-fold tensor product
of a given vertex operator algebra with itself and a permutation
automorphism of the tensor factors.  The latter construction is based
on a certain change of variables in the worldsheet geometry in the
sense of string theory.  In the case of a lattice that is the
orthogonal direct sum of copies of a given lattice, these two very
different constructions can both be carried out, and must produce
isomorphic twisted modules, by a theorem of the first author jointly
with Dong and Mason.  In this paper, we explicitly construct an
isomorphism, thereby providing, {}from both mathematical and physical
points of view, a direct link between space-time geometry and
worldsheet geometry in this setting.
\end{abstract}

\maketitle

\section{Introduction and preliminaries}

Twisted modules for vertex operator algebras arose in the work of the
third author with I. Frenkel and A. Meurman \cite{FLM1}, \cite{FLM2},
\cite{FLM3} for the case of a lattice vertex operator algebra and the
lattice isometry $-1$, in the course of the construction of the
moonshine module vertex operator algebra. This structure came to be
understood as an ``orbifold model" in the sense of conformal field
theory and string theory.  Twisted modules are the mathematical
counterpart of ``twisted sectors", which are the basic building blocks
of orbifold models in conformal field theory and string theory (see
\cite{DHVW1}, \cite{DHVW2}, \cite{DFMS}, \cite{DVVV}, \cite{DGM}, as
well as \cite{KS}, \cite{FKS}, \cite{Ba1}, \cite{Ba2}, \cite{BHS},
\cite{dBHO}, \cite{HO}, \cite{GHHO}, \cite{Ba3} and \cite{HH}).  The
notion of twisted module for a vertex operator algebra is a
generalization of the notion of module in which the action of an
automorphism of the vertex operator algebra is incorporated.  Given a
vertex operator algebra and an automorphism, it is an open problem as
to how to construct a corresponding twisted module in general.
However, the problem of constructing twisted modules has been solved
in particular for two families of vertex operators and their
automorphisms.  One of these constructions is that initiated by the
third author \cite{L1} in the case of a lattice vertex operator
algebra and an automorphism arising {}from an arbitrary lattice
isometry, generalizing the joint work of the third author mentioned
above.  This construction is ultimately based on the lattice isometry,
and thus, {}from a physical point of view, is related to the
space-time geometry associated with the lattice in the sense of string
theory.  The other construction is due to the first author, jointly
with C. Dong and G. Mason \cite{BDM}, in the case of a multi-fold
tensor product of a given vertex operator algebra with itself and a
permutation automorphism of the tensor factors.  The latter
construction is based on a change of variables in the worldsheet
geometry in the sense of string theory.  Now, in the case of a lattice
which is the orthogonal direct sum of copies of a given lattice, these
two very different constructions can both be carried out. By a theorem
of the first author jointly with Dong and Mason \cite{BDM}, in this
case these two constructions must produce isomorphic twisted modules.
In this paper, we explicitly construct an isomorphism, thereby
providing, {}from both mathematical and physical points of view, a
direct link between space-time geometry and worldsheet geometry in
this interesting setting.

The precise notion of vertex operator algebra was developed in
\cite{FLM3}, following Borcherds' introduction of the notion of vertex
algebra in \cite{Bo}.  Twisted vertex operators were discovered and
used in \cite{LW}.  The first orbifold conformal field theory (as it
came to be understood) was introduced in \cite{FLM1}.  Formal calculus
arising {}from twisted vertex operators associated to an even lattice
was systematically developed in \cite{L1}, \cite{FLM2}, \cite{FLM3}
and \cite{L2}, and the twisted Jacobi identity was formulated and
shown to hold for these operators (see also \cite{DL}).  These results
led to the introduction of the notion of $g$-twisted $V$-module
\cite{FFR}, \cite{D2}, for $V$ a vertex operator algebra and $g$ an
automorphism of $V$.  This notion records the properties of twisted
operators obtained in \cite{L1}, \cite{FLM1}, \cite{FLM2}, \cite{FLM3}
and \cite{L2}, and provides an axiomatic definition of the notion of
twisted sectors.  In general, given a vertex operator algebra $V$ and
an automorphism $g$ of $V$, it is an open problem as to how to
construct a $g$-twisted $V$-module.

In \cite{BDM}, twisted modules for a permutation acting on a tensor
product vertex operator algebra were constructed and classified.  Let
$V$ be a vertex operator algebra, and for a fixed positive integer
$k$, consider the tensor product vertex operator algebra $V^{\otimes
k}$ (cf. \cite{FHL}).  Any element of the symmetric group acts on
$V^{\otimes k}$ in the obvious way, and this is the setting for
permutation-twisted modules.  {}From the physical point of view, this
is the setting for permutation orbifold theory and has been studied,
for example, in \cite{KS}, \cite{FKS}, \cite{Ba1}, \cite{Ba2},
\cite{BHS} and \cite{Ba3}.  In the case of $V$ a lattice vertex
operator algebra, the construction of \cite{BDM} becomes a special
case of the more general results of \cite{L1}, \cite{FLM2}, \cite{L2},
and \cite{DL}, and this overlap of constructions is the basis for this
paper.

In this paper, we begin by giving some preliminary definitions,
including the axiomatic definition of twisted module, and review the
construction of a vertex operator algebra $V_L$ associated to a
positive-definite even lattice $L$. In Section \ref{specialize-section},
we give the setting for permutation-twisted modules associated to a
lattice that is the orthogonal direct sum of copies of a positive-definite
even lattice.  For $K$ a positive-definite even lattice, $k$ a positive
integer, and $L = K \oplus K \oplus \cdots \oplus K$ the orthogonal direct
sum of $k$ copies of $K$, we consider the lattice automorphism of $L$ given
by permuting the direct sum factors $K$ by a cyclic permutation $\nu =
(1 \; 2\; \cdots k)$.  The lattice automorphism $\nu$ lifts to the
$V_L$-automorphism $\hat{\nu}$ given by cyclicly permuting the $k$
tensor copies of $V_K$ in $V_L = V_K^{\otimes k}$.

In Section \ref{spacetime-section}, we give the construction of
irreducible $\hat{\nu}$-twisted $V_L$-modules in this setting,
following \cite{L1}, \cite{FLM2}, \cite{DL}, and we calculate
the graded dimensions of these modules.  In Section
\ref{bdm-section}, we present the construction of
$\hat{\nu}$-twisted $V_L$-modules following \cite{BDM} specialized to
this setting.  We also recall the results of \cite{BDM} pertaining to
the determination of the irreducible $\hat{\nu}$-twisted
$V_L$-modules.  In particular, we note that in \cite{BDM} it is shown
that the category of irreducible $\hat{\nu}$-twisted $V_L$-modules is
isomorphic to the category of irreducible $V_K$-modules.  In Section
\ref{compare-section}, we use this determination of
$\hat{\nu}$-twisted $V_L$-modules to conclude that the
$\hat{\nu}$-twisted $V_L$-module constructed via the method of
\cite{L1}, \cite{FLM2} must be isomorphic to some $\hat{\nu}$-twisted
$V_L$-module constructed via the method of \cite{BDM}.  We then recall
the classification of irreducible $V_K$-modules given in \cite{D1}
(see also \cite{DLiM1}; cf. \cite{LL}).  We use this to prove that
under the isomorphism of categories {}from \cite{BDM} the irreducible
$V_K$-module corresponding to the $\hat{\nu}$-twisted $V_L$-module
following the \cite{L1} construction must be $V_K$ itself.  We prove
this using graded dimensions.  This allows us to pick out which
$\hat{\nu}$-twisted $V_L$-module under the \cite{BDM} construction
must be isomorphic to the $\hat{\nu}$-twisted $V_L$-module of the
\cite{L1} construction.  In Section \ref{explicit-section}, using the
existence of the isomorphism between the two constructions of
$\hat{\nu}$-twisted $V_L$-modules proved in Section
\ref{compare-section}, we explicitly determine the isomorphism.  We
also show how to generalize these results to $g$-twisted
$V_L$-modules, where $g$ is an arbitrary permutation on $k$ letters,
and in addition, to arbitrary irreducible modules and twisted modules,
corresponding to cosets of the relevant lattices.

We would now like to comment on some motivations and expected
implications of this work.  As mentioned previously, the construction
of $\hat{\nu}$-twisted $V_L$-modules following \cite{L1}, \cite{FLM1},
\cite{FLM2}, \cite{FLM3} and \cite{DL} is ultimately based on the
lattice automorphism $\nu$ and thus, {}from the physical point of
view, is inherently based on the ``orbifolding" of the space-time
geometry in the sense of string theory.  In fact, the lattice vertex
operator algebra $V_{L}$ and $V_{L}$-modules can be interpreted
physically by quantizing the classical theory of strings propagating
in the space-time torus $\mathbb{R}^{{\rm rank}\; L}/L$.  In
this picture, strings in the torus are studied as strings in the
Euclidean space satisfying periodic boundary conditions.  Twisted
modules for $V_{L}$ can be analogously interpreted physically by
quantizing the classical theory of strings propagating in the orbifold
obtained by taking the quotient of the torus by a group.  Strings in
this orbifold can be studied as strings in the Euclidean space
satisfying ``periodic boundary conditions up to actions of elements of
the group."  However, the twisted modules are quite subtle to
construct mathematically.  The mathematical construction of twisted
modules for $V_{L}$ in \cite{L1}, \cite{FLM1}, \cite{FLM2},
\cite{FLM3} and \cite{DL} can in fact be physically interpreted using
this space-time picture; indeed, see \cite{DHVW1} and the related
string-theoretic works on strings on orbifolds, and on orbifold models
in conformal field theory.

On the other hand, the construction of $\hat{\nu}$-twisted
$V_L$-modules following \cite{BDM} is, in general, independent of the
given lattice and instead relies on an operator derived {}from a
change of coordinates related to the conformal geometry of propagating
strings, and thus, {}from the physical point of view, is based on the
the worldsheet geometry; see Remark \ref{worldsheet-remark}.  In fact,
the ``periodic boundary conditions up to actions of elements of the
group" mentioned above shows that one needs to consider multivalued
functions on the worldsheet of strings in an orbifold.  Such
multivalued functions are exactly the ones used in \cite{BDM} to
construct twisted modules in the setting of that work.

The completely different geometric foundations for the two
constructions highlight just how different these two constructions
are, and yet they give isomorphic $\hat{\nu}$-twisted $V_L$-modules.
Thus, {}from both mathematical and physical points of view, the
isomorphism between the two constructions gives a direct link between
space-time geometry and worldsheet geometry in this interesting
setting.  {}From a purely mathematical viewpoint, one of the important
applications and also one of the motivations for this work is that
this isomorphism between the ``space-time" construction (\cite{L1},
\cite{FLM2}, \cite{DL}) and the ``worldsheet" construction \cite{BDM}
provides a way of transporting interesting structures that have been
developed following the ``space-time" construction to the conformal
geometry of the worldsheet.  For instance, we expect to relate the
present work to the work of the third author jointly with Doyon and
Milas in \cite{DLeM1} and \cite{DLeM2}.

{\noindent} {\bf Acknowledgments} The first author would like to thank
the University of Notre Dame for the research leave that facilitated
this collaboration, and would also like to thank Rutgers University,
and Bill and Kathy Exner for their hospitality.  The authors thank the
referee for helpful comments.  The authors gratefully acknowledge
partial support {}from NSF grant DMS-0401302.

{\noindent} {\bf Notation} $\mathbb{Z}_+$ denotes the positive
integers and $\mathbb{N}$ denotes the nonnegative integers.

\subsection{Vertex operator algebras, modules, automorphisms and twisted
modules}\label{generalities-section}

In this section, we review the definitions of vertex operator algebra
and $g$-twisted $V$-module for a vertex operator algebra $V$ and an
automorphism $g$ of $V$ of finite order.

We begin by recalling the notion of vertex operator algebra, following
the notation and terminology of \cite{FLM2} and \cite{LL}. Let $x,
x_0, x_1, x_2,$ etc., denote commuting independent formal variables.
Let $\delta (x) = \sum_{n \in \Z} x^n$.  We will use the binomial
expansion convention, namely, that any expression such as $(x_1 -
x_2)^n$ for $n \in \C$ is to be expanded as a formal power series in
nonnegative integral powers of the second variable, in this case
$x_2$.

A {\it vertex operator algebra} is a $\Z$-graded vector space
\begin{equation}
V=\coprod_{n\in\Z}V_n
\end{equation}
satisfying ${\rm dim} \, V < \infty$ and $V_n = 0$ for $n$ sufficiently negative
and equipped with a linear map
\begin{eqnarray}
V &\longrightarrow& (\mbox{End}\,V)[[x,x^{-1}]]\\
v &\mapsto& Y(v,x)= \sum_{n\in\Z}v_nx^{-n-1} \nonumber
\end{eqnarray}
and with two distinguished vectors ${\bf 1}\in V_0$, (the {\it vacuum vector})
and $\omega\in V_2$
(the {\it conformal element}) satisfying the following conditions for $u, v \in
V$:
\begin{eqnarray}
& &u_nv=0\ \ \ \ \ \mbox{for $n$ sufficiently large};\\
& &Y({\bf 1},x)=1;\\
& &Y(v,x){\bf 1}\in V[[x]]\ \ \ \mbox{and}\ \ \ \lim_{x\to
0}Y(v,x){\bf 1}=v;
\end{eqnarray}
\begin{eqnarray}
\lefteqn{x^{-1}_0\delta\left(\frac{x_1-x_2}{x_0}\right) Y(u,x_1)Y(v,x_2) -
x^{-1}_0\delta\left(\frac{x_2-x_1}{-x_0}\right) Y(v,x_2)Y(u,x_1)}\\
&=& x_2^{-1}\delta \left(\frac{x_1-x_0}{x_2}\right) Y(Y(u,x_0)v,x_2)
\hspace{1.3in} \nonumber
\end{eqnarray}
(the {\it Jacobi identity});
\begin{equation}
[L(m),L(n)]=(m-n)L(m+n)+\frac{1}{12}(m^3-m)\delta_{m+n,0}c
\end{equation}
for $m, n\in \Z,$ where
\begin{equation}
L(n)=\omega_{n+1}\ \ \ \mbox{for $n\in \Z$, \ \ \ \ i.e.},\
Y(\omega,x)=\sum_{n\in\Z}L(n)x^{-n-2}
\end{equation}
and $c \in \mathbb{C}$ (the {\it central charge} of $V$);
\begin{eqnarray}
& &L(0)v=nv=(\mbox{wt}\,v)v \ \ \ \mbox{for $n \in \Z$ and $v\in V_n$}; \\
& &\frac{d}{dx}Y(v,x)=Y(L(-1)v,x).
\end{eqnarray}
This completes the definition. We denote the vertex operator algebra just
defined by $(V,Y,{\bf 1},\omega)$ (or briefly, by $V$).

The {\it graded dimension} of a vertex operator algebra $V = \coprod_{n \in \Z}
V_n$ is defined to be
\begin{equation}
{\rm dim}_* V = {\rm tr}_V q^{L(0) - c/24} = q^{-c/24} \sum_{n \in \Z} ({\rm
dim} V_n) q^n
\end{equation}
where $q$ is a formal variable and $c$ is the central charge of $V$.

An {\it automorphism} of a vertex operator algebra $V$ is a linear
automorphism $g$ of $V$ preserving ${\bf 1}$ and $\omega$ such that
the actions of $g$ and $Y(v,x)$ on $V$ are compatible in the sense
that
\begin{equation}
g Y(v,x) g^{-1}=Y(gv,x)
\end{equation}
for $v\in V.$ Then $g V_n\subset V_n$ for $n\in \mathbb{Z}$.  If $g$
has finite order, $V$ is a direct sum of the eigenspaces $V^j$ of $g$,
\begin{equation}
V=\coprod_{j\in \Z /k \Z }V^j,
\end{equation}
where $k \in \mathbb{Z}_+$ is a period of $g$ (i.e., $g^k = 1$ but $k$
is not necessarily the order of $g$) and
\begin{equation}
V^j=\{v\in V \; | \; g v= \eta^j v\},
\end{equation}
for $\eta$ a fixed primitive $k$-th root of unity.

We next recall the notion of $g$-twisted $V$-module, which records the
properties of twisted vertex operators obtained in \cite{L1},
\cite{FLM2} and \cite{L2}.  We follow the notation and terminology of
\cite{BDM}.  Let $(V,Y,{\bf 1},\omega)$ be a vertex operator algebra
and let $g$ be an automorphism of $V$ of period $k \in
\mathbb{Z}_+$. A {\it $g$-twisted $V$-module} $M$ is a $\C$-graded
vector space
\begin{equation}
M=\coprod_{\lambda \in \C}M_\lambda
\end{equation}
such that for each $\lambda$, $\dim M_{\lambda}< \infty $ and $M_{n/k
+\lambda}=0$ for all sufficiently negative integers $n$.  In addition,
$M$ is equipped with a linear map
\begin{eqnarray}
V &\longrightarrow & (\mbox{End}\,M)[[x^{1/k},x^{-1/k}]] \\
v &\mapsto & Y^g (v,x)=\sum_{n\in \frac{1}{k}\Z }v_n^g x^{-n-1} \nonumber
\end{eqnarray}
satisfying the following conditions for $u,v\in V$ and $w\in M$:
\begin{eqnarray}
\label{cosetrelation}
& &Y^g(v,x)=\sum_{n\in j/k+ \Z}v_n^g x^{-n-1}\ \ \ \
\mbox{for $j\in \Z/k\Z$ and $v\in V^j$}; \\
& &v_n^g w=0\ \ \ \mbox{for $n$ sufficiently large};\\
& &Y^g ({\bf 1},x)=1;
\end{eqnarray}
\begin{multline}
x^{-1}_0\delta\left(\frac{x_1-x_2}{x_0}\right)
Y^g(u,x_1)Y^g(v,x_2)-x^{-1}_0\delta\left(\frac{x_2-x_1}{-x_0}\right) Y^g
(v,x_2)Y^g (u,x_1)\\
= x_2^{-1}\frac{1}{k}\sum_{j\in \Z /k \Z}
\delta\left(\eta^j\frac{(x_1-x_0)^{1/k}}{x_2^{1/k}}\right)Y^g (Y(g^j
u,x_0)v,x_2)
\end{multline}
(the {\it twisted Jacobi identity}) where $\eta$ is a fixed primitive
$k$-th root of unity;
\begin{equation}\label{Viralgrelations}
[L^g(m),L^g(n)]=(m-n)L^g (m+n)+\frac{1}{12}(m^3-m)\delta_{m+n,0}c
\end{equation}
for $m, n\in \Z$, where $c$ is the central charge of $V$, and
\begin{equation}\label{define-L^g}
L^g (n)=\omega_{n+1}^g \ \ \ \mbox{for $n\in \Z$, \ \ \ \ \ i.e.},\ Y^g
(\omega,x)=\sum_{n\in \Z} L^g (n)x^{-n-2};
\end{equation}
\begin{equation}\label{L^g-grading}
L^g (0) w=\lambda w \qquad \mbox{for $w \in M_\lambda$};
\end{equation}
\begin{equation}
\frac{d}{dx}Y^g (v,x)=Y^g (L(-1)v,x).
\end{equation}
(Formula (\ref{cosetrelation}) can be expressed as follows: For $v \in
V$,
\begin{equation}
Y^g (gv,x) = \lim_{x^{1/k} \rightarrow \eta^{-1} x^{1/k}} Y^g(v, x),
\end{equation}
where the limit stands for formal substitution.)  This completes the
definition of $g$-twisted $V$-module.  We denote such a module by
$(M,Y^g)$ (or briefly, by $M$).

If we take $g=1$, then we obtain the notion of (ordinary) $V$-module.

We call a $g$-twisted $V$-module $M$ {\it simple} or {\it irreducible}
if the only submodules are 0 and $M$.

A vertex operator algebra is {\it simple} if it is simple as a module
for itself.

Note that the notion of graded dimension still makes sense for
$g$-twisted $V$-modules (and thus for ordinary $V$-modules); that is,
we have
\begin{equation*}
{\rm dim}_* M = {\rm tr}_M q^{L^g(0) - c/24}.
\end{equation*}

Let $(M^1, Y_1^g)$ and $(M^2, Y_2^g)$ be $g$-twisted $V$-modules.  A
{\it homomorphism} {}from $M^1$ to $M^2$ is a linear map $f: M^1
\longrightarrow M^2$ such that
\begin{equation}
f(Y_1 ^g(v,x) w) = Y_2^g(v,x) f(w)
\end{equation}
for $v \in V$ and $w \in M^1$.

\subsection{Lattice vertex operator algebras}\label{lattice-section}

We next recall the construction of vertex operator algebras and
related structures corresponding to a lattice equipped with an
isometry, following the notation and terminology of \cite{FLM3} and
using the setting and results of \cite{L1} and \cite{FLM2}.

Let $L$ be a positive-definite even lattice, with (nondegenerate
symmetric) $\mathbb{Z}$-bilinear form $\langle \cdot, \cdot \rangle$.
(There should be no confusion between this use of the symbol $L$ and
the operators $L(n)$.)  Let $\nu$ be an isometry of $L$, and let $k
\in \mathbb{Z}_+$ such that the following hold:
\begin{equation}
\nu^k = 1;
\end{equation}
if $k$ is even, then
\begin{equation}\label{even-equation}
\langle \nu^{k/2} \alpha, \alpha \rangle \in 2 \mathbb{Z} \ \ \ \ \mbox{for
$\alpha \in L$}
\end{equation}
(which can be arranged by doubling $k$ if necessary).
Observe that under these assumptions,
\begin{equation}\label{even-condition}
\left< \sum_{j = 0}^{k-1} \nu^j \alpha, \alpha \right> \in 2 \mathbb{Z}
\end{equation}
for $\alpha \in L$.  Let $\eta$ be a fixed primitive $k$-th root of unity.
Define the functions
\begin{eqnarray}\label{commutator-definition-0}
C_0: L \times L &\longrightarrow& \C^\times \\
(\alpha, \beta) &\mapsto& (-1)^{ \langle \alpha,  \beta \rangle} \nonumber,
\end{eqnarray}
and
\begin{eqnarray}\label{commutator-definition}
C: L \times L &\longrightarrow& \C^\times \\
(\alpha, \beta) &\mapsto& (-1)^{\sum_{j = 0}^{k-1} \langle \nu^j \alpha,  \beta
\rangle} \eta^{\sum_{j = 0}^{k-1} \langle j \nu^j \alpha , \beta \rangle}
\nonumber \\
& &\quad = \prod_{j=0}^{k-1} (-\eta^{j} )^{\langle \nu^j \alpha, \beta
\rangle}. \nonumber
\end{eqnarray}
Note that $C_0$ and $C$ are bilinear into the abelian group $\C^\times$; i.e.,
\begin{eqnarray*}
C(\alpha + \beta, \gamma) = C(\alpha,\gamma)C(\beta, \gamma) \\
C(\alpha, \beta + \gamma) = C(\alpha, \beta)C(\alpha, \gamma)
\end{eqnarray*}
for $\alpha, \beta, \gamma \in L$, and similarly for $C_0$.  By the fact that
$L$ is even, we have $C_0 (\alpha, \alpha) = 1$, and by (\ref{even-condition}),
we have $C(\alpha, \alpha) = 1$.  We also note that both $C_0$ and $C$ are
$\nu$-invariant, that is, $C(\nu \alpha, \nu \beta) = C(\alpha, \beta)$ and
similarly for $C_0$.  Moreover, $C(\beta, \alpha) = C(\alpha, \beta)^{-1}$.

Set
\begin{equation}
\eta_0 = (-1)^k \eta.
\end{equation}
Then $\eta_0$ is a primitive $2k$-th root of unity if $k$ is odd, and
$-1$ and $\eta$ are powers of $\eta_0$ for any $k$.

The maps $C_0$ and $C$ determine uniquely (up to equivalence) two central
extensions of $L$ by the cyclic group $\langle \eta_0 \rangle$,
\begin{equation}\label{exact-0}
1 \rightarrow \langle \eta_0 \rangle \rightarrow \hat{L} \bar{\longrightarrow} L
\rightarrow 1,
\end{equation}
\begin{equation}\label{exact}
1 \rightarrow \langle \eta_0 \rangle \rightarrow \Lnu \bar{\longrightarrow} L
\rightarrow 1,
\end{equation}
with commutator maps $C_0$ and $C$, respectively,
i.e.,
such that
\begin{eqnarray}\label{commutator=C0}
aba^{-1} b^{-1} &=& C_0( \bar{a}, \bar{b})  \qquad \mathrm{for} \quad a,b \in
\hat{L} ,\\
aba^{-1} b^{-1} &=& C( \bar{a}, \bar{b})  \qquad \mathrm{for} \quad a,b \in
\Lnu .
\end{eqnarray}
There is a natural set-theoretic identification (which is not an
isomorphism of groups unless $k = 1$ or $k = 2$) between the groups
$\hat{L}$ and $\Lnu$ such that the respective group multiplications
$\times$ and $\times_\nu$ are related by
\begin{equation}\label{identify-central-extensions}
a \times b = \prod_{0<j<k/2} (- \eta^j)^{\langle \nu^{-j} \bar{a}, \bar{b}
\rangle} a \times_\nu b \qquad \mathrm{for} \quad a,b \in \hat{L}.
\end{equation}

Observe further that since $C_0$ is $\nu$-invariant, if we replace the
map $\ \bar{} \ $ in (\ref{exact-0}) by $\nu \circ \ \bar{} \ $, we
obtain another central extension of $L$ by $\langle \eta_0 \rangle$
with commutator map $C_0$.  By uniqueness of the central extension of
$L$, there is an automorphism $\hat{\nu}$ of $\hat{L}$ (fixing
$\eta_0$) such that $\hat{\nu}$ is a lifting of $\nu$, i.e., such that
\begin{equation}\label{nu-on-hatL-0}
(\hat{\nu} a)^{\bar{}} = \nu \bar{a} \quad \mathrm{for} \quad a \in \hat{L} .
\end{equation}
The map $\hat{\nu}$ is also an automorphism of $\Lnu$ satisfying
\begin{equation}\label{nu-on-hatL}
(\hat{\nu} a)^{\bar{}} = \nu \bar{a} \quad \mathrm{for} \quad a \in \Lnu .
\end{equation}
Moreover, we may choose the lifting $\hat{\nu}$ of $\nu$ so that
\begin{equation}\label{nuhata=a}
\hat{\nu} a = a \quad \mathrm{if} \quad \nu \bar{a} = \bar{a}
\end{equation}
(see (\ref{propertyofnuhat}) below), and we have
\begin{equation}\label{nuhat^k=1}
\hat{\nu}^k = 1,
\end{equation}
a nontrivial fact (see \cite{L1}).

We now use the central extension $\hat{L}$ to construct a vertex
operator algebra $V_L$ equipped with an automorphism $\hat{\nu}$ of
period $k$, induced {}from the automorphism $\hat{\nu}$ of $\hat{L}$.
In Section \ref{specialize-section} we will specialize our setting, specifying
$\nu$ and
$\hat{\nu}$ in this setting.  Then in Section \ref{spacetime-section}
we will use the central extension $\Lnu$ in the specialized setting to
construct an irreducible $\hat{\nu}$-twisted module for the vertex operator
algebra
$V_L$, following \cite{L1}, \cite{FLM2}, \cite{DL} implemented in our
specialized setting.  (In \cite{L1}, \cite{FLM2}, \cite{DL}, such
$\hat{\nu}$-twisted modules are constructed in the general case.)

Embed $L$ canonically in the $\C$-vector space $\h = \C
\otimes_\mathbb{Z} L$, and extend the $\Z$-bilinear form on $L$ to a
$\C$-bilinear form $\langle \cdot, \cdot \rangle$ on $\h$.  The
corresponding affine Lie algebra is
\begin{equation}
\hat{\h} = \h \otimes \C[t,t^{-1}] \oplus \C {\bf c},
\end{equation}
with brackets determined by
\begin{equation}
[\alpha \otimes t^m, \beta \otimes t^n]=\langle \alpha , \beta \rangle
m\delta_{m+n,0}{\bf c}\ \ {\rm for}\ \ \alpha, \beta \in{\h}, \ \ m,n\in \Z,
\end{equation}
\begin{equation}
[{\bf c},\hat{\h}]=0.
\end{equation}
Then $\hat{\h}$ has a $\Z$-gradation, the {\it weight gradation},
given by
\begin{equation}
{\rm wt}\,(\alpha \otimes t^n)=-n \ \ {\rm and} \ \ {\rm wt}\, {\bf c}=0
\end{equation}
for $\alpha \in {\h}$ and $n\in \Z$.

Set
\begin{equation}
\hat{\h}^+={\h}\otimes t \C[t] \ \ {\rm and} \ \ \hat{\h}^-={\h}\otimes t^{-1}
\C[t^{-1}].
\end{equation}
The subalgebra
\begin{equation}
\hat{\h}_{\Z} = \hat{\h}^+\oplus\hat{\h}^-\oplus \C {\bf c}
\end{equation}
of $\hat{\h}$ is a Heisenberg algebra, in the sense that its
commutator subalgebra equals its center, which is one-dimensional.
Consider the induced $\hat{\h}$-module, irreducible even under
$\hat{\h}_{\Z}$, given by
\begin{equation}
M(1)=U(\hat{\h})\otimes_{U( \h \otimes \C [t] \oplus \C {\bf c})} \C \simeq
S(\hat{\h}^-)\ \ \ (\mbox{linearly}),
\end{equation}
where ${\h}\otimes \C [t]$ acts trivially on $\C$ and ${\bf c}$ acts as 1,
$U(\cdot)$ denotes universal enveloping algebra and $S(\cdot)$ denotes
symmetric algebra.  The $\hat{\h}$-module $M(1)$ is $\Z$-graded so
that wt\,1\ =\ 0 (we write 1 for $1\otimes 1$):
\begin{equation}
M(1)=\coprod_{n\in \mathbb{N}} M(1)_n,
\end{equation}
where $M(1)_n$ denotes the homogeneous subspace of weight $n$.

Form the induced $\hat{L}$-module and $\C$-algebra
\begin{eqnarray}
\C \{L\} &=& \C [\hat{L}] \otimes_{\mathbb{C}[ \langle \eta_0 \rangle] } \C\\
&\simeq &\C [L]\ \ \ \qquad \qquad (\mbox{linearly}), \nonumber
\end{eqnarray}
where $\C [\cdot]$ denotes group algebra.  For $a \in \hat{L}$, write
$\iota(a)$ for the image of $a$ in $\C \{L\}$. Then the action of $\hat{L}$
on $\C\{L\}$ is given by
\begin{equation}
a\cdot\iota(b)=\iota(a)\iota(b)=\iota(ab)
\end{equation}
for $a,b\in\hat{L}$. We give $\C \{L\}$ the $\C$-gradation determined by:
\begin{equation}
\mbox{wt}\,\iota(a)=\frac{1}{2}\langle \bar{a},\bar{a}\rangle \ \ \ \
\mbox{for}\ \
a\in \hat{L}.
\end{equation}
Also define a grading-preserving action of ${\h}$ on $\C\{L\}$ by:
\begin{equation}
h\cdot\iota(a)= \langle h,\bar{a}\rangle \iota(a)
\end{equation}
for $h\in{\h}$, and define
\begin{equation}
x^h\cdot\iota(a) = x^{\langle h,\bar{a}\rangle }\iota(a)
\end{equation}
for $h\in{\h}.$

Set
\begin{eqnarray}
V_L &=& M(1)\otimes_{\C} \C \{L\}\\
&\simeq & S(\hat{\h}^-)\otimes \C [L] \ \ \ \ \qquad (\mbox{linearly}) \nonumber
\end{eqnarray}
and give $V_L$ the tensor product $\C$-gradation:
\begin{equation}
V_L=\coprod_{n\in \C}(V_L)_n.
\end{equation}
We have wt$\,\iota(1)=0,$ where we identify $\C \{L\}$ with $1\otimes
\C \{L\}$.  Then $\hat{L}$, $\hat{\h}_{\Z}$, $h$, $x^h$ $(h\in{\h})$
act naturally on $V_L$ by acting on either $M(1)$ or $\C \{L\}$ as
indicated above. In particular, ${\bf c}$ acts as 1.

For $\alpha \in \h$, $n\in \Z$, we write $\alpha(n)$ for the operator
on $V_L$ determined by $\alpha\otimes t^n$. For $\alpha\in \h,$ set
\begin{equation}
\alpha(x)=\sum_{n\in \Z} \alpha(n) x^{-n-1}.
\end{equation}
We use a normal ordering procedure, indicated by open colons, which
signify that the enclosed expression is to be reordered if necessary
so that all the operators $\alpha(n)$ ($\alpha\in \h$, $n<0$) and
$a\in \hat{L}$ are to be placed to the left of all the operators
$\alpha(n),$ and $x^{\alpha}$ ($\alpha \in \h$, $n\ge 0$) before the
expression is evaluated. For $a \in \hat{L}$, set
\begin{equation}
Y(a,x)= \ _\circ^\circ \ e^{\int(\bar{a}(x)- \bar{a}(0)x^{-1})} a x^{\bar{a}}
\ _\circ^\circ ,
\end{equation}
using an obvious formal integration notation.  Let $a \in \hat{L}$,
$\alpha_1, \dots, \alpha_m \in \h$, $n_1,\dots, n_m \in \Z_+$ and set
\begin{eqnarray}
v &=& \alpha_1(-n_1)\cdots \alpha_m(-n_m) \otimes \iota(a)\\
&=& \alpha_1(-n_1)\cdots \alpha_m(-n_m) \cdot \iota(a) \in V_L. \nonumber
\end{eqnarray}
Define
\begin{multline}
Y(v,x) = \ _\circ^\circ  \left(\frac{1}{(n_1-1)!}\left(\frac{
d}{dx}\right)^{n_1-1} \alpha_1(x)\right)\cdots
\\\left(\frac{1}{(n_m-1)!}\left(\frac{d}{dx}\right)^{n_m-1}\alpha_m(x)\right)Y(a,x)
\ _\circ^\circ .
\end{multline}
This gives us a well-defined linear map
\begin{eqnarray}
V_L&\rightarrow&(\mbox{End}\,V_L)[[x, x^{-1}]] \\
v &\mapsto& Y(v,x)=\displaystyle{\sum_{n\in \Z}}v_nx^{-n-1}, \ \ \
v_n\in\mbox{End}\,V_L. \nonumber
\end{eqnarray}
Set ${\bf 1} = 1 = 1 \otimes 1 \in V_L$ and
\begin{equation}\label{define-omega}
\omega = \frac{1}{2} \sum_{i=1}^{{\rm dim} \, {\h}} h_i(-1)h_i(-1) {\bf 1},
\end{equation}
where $\{h_i\}$ is an orthonormal basis of $\h$. Then $V_L = (V_L, Y,
{\bf 1}, \omega)$ is a simple vertex operator algebra of central charge
$$c = {\rm dim} \, {\h} = {\rm rank} \, L.$$

\begin{rema}\label{independenceofchoices}
{\em The construction of the vertex operator algebra $V_L$ depends on
the central extension (\ref{exact-0}) subject to
(\ref{commutator=C0}), and hence on the choices of $k\in \Z_+$ and the
primitive root of unity $\eta$.  But it is a standard fact that $V_L$
is independent of these choices, up to isomorphism of vertex operator
algebras preserving the $\hat{\h}$-module structure; see Proposition
6.5.5, and also Remarks 6.5.4 and 6.5.6, of \cite{LL}.  In particular,
$V_L$ as constructed above is essentially the same as $V_L$
constructed {}from a central extension of the type (\ref{exact-0})
subject to (\ref{commutator=C0}) but with the kernel of the central
extension replaced by the group $\langle \pm 1 \rangle$.  For the
purpose of constructing twisted modules, it is valuable to have this
flexibility.  We will use these properties of lattice vertex operator
algebras below.}
\end{rema}

Following  \cite{L1} (and see also \cite{DL}), we note that the
automorphism $\nu$ of $L$ acts in a natural way on ${\h},$ on
$\hat{\h}$ (fixing $c$) and on $M(1)$, preserving the gradations, and
for $u\in \hat{\h}$ and $m\in M(1),$
\begin{equation}
\nu (u\cdot m)=\nu (u)\cdot \nu(m).
\end{equation}
The automorphism $\nu$ of $L$ lifted to the automorphism $\hat{\nu}$ of
$\hat{L}$
satisfies
\begin{equation}
\hat{\nu}(h\cdot \iota(a))=\nu(h)\cdot \hat{\nu}\iota(a),
\end{equation}
for $h \in \h$ and $a \in \hat{L}$, and we have
\begin{equation}
\hat{\nu}(\iota(a)\iota(b))=\hat{\nu} (a\cdot \iota(b)) = \hat{\nu}(a) \cdot
\hat{\nu} \iota(b) = \hat{\nu} \iota(a) \hat{\nu} \iota(b),
\end{equation}
\begin{equation}
\hat{\nu}(x^h \cdot \iota(a)) = x^{\nu(h)} \cdot \hat{\nu}\iota(a).
\end{equation}
Thus we have a natural grading-preserving automorphism of $V_L$, which
we also call $\hat{\nu}$, which acts via $\nu\otimes\hat{\nu},$ and this action
is compatible with the other actions:
\begin{eqnarray}
\hat{\nu}(a\cdot v) &=& \hat{\nu}(a)\cdot\hat{\nu}(v)\\
\hat{\nu}(u\cdot v) &=& \nu(u)\cdot \hat{\nu}(v) \\
\hat{\nu}(x^h\cdot v) &=& x^{\nu(h)}\cdot \hat{\nu}(v)
\end{eqnarray}
for $a\in \hat{L}$, $u\in \hat{\h}$, $h\in {\h}$, and $v\in V_L$, so
that $\hat{\nu}$ is an automorphism of the vertex operator algebra
$V_L$.  In Section \ref{specialize-section}, we will specialize this
general setting to a specific type of lattice $L$ and automorphism
$\nu$ and use the automorphism $\hat{\nu}$ of $V_L$ to construct a
$\hat{\nu}$-twisted $V_L$-module.

Finally, in the general setting, we consider the graded dimension of
the vertex operator algebra $V_L$.  Let $\Theta_L(q)$ be the theta
function corresponding to $L$; that is,
\begin{equation}
\Theta_L(q) = \sum_{\alpha \in L} q^{ \langle \alpha , \alpha \rangle /2},
\end{equation}
and let $\eta(q)$ be the Dedekind eta function, given by
\begin{equation}
\eta(q) = q^{1/24} \prod_{n \in \Z_+} (1 - q^n) .
\end{equation}
Then we have
\begin{equation}
{\rm dim}_* V_L = \frac{\Theta_L(q)}{\eta(q)^d} .
\end{equation}

\section{Specialization of the general setting and the "space-time"
construction of  $\hat{\nu}$-twisted $V_L$-modules}\label{specialize-section}

In \cite{L1}, twisted modules for the vertex operator algebra
associated to a positive definite lattice and a lattice isometry are
constructed.  In \cite{BDM}, twisted modules for a vertex operator
algebra which is the $k$-fold tensor product ($k \in \Z_+$) of a
vertex operator algebra $V$ with itself, twisted by a permutation
automorphism of $V^{\otimes k}$, are constructed.  In this paper, we
investigate these two constructions in a setting in which they
overlap.  We will now describe this setting, specializing the general
notions above.  In Section \ref{spacetime-section}, we will carry out
the construction of \cite{L1} in this setting, and in Section
\ref{bdm-section} we will carry out the construction of \cite{BDM} in
this setting.

\subsection{The setting}
Let $K$ be a positive definite even lattice with symmetric
bilinear form given by $\langle \cdot, \cdot \rangle$, and for a fixed
$k \in \Z_+$, let
\begin{equation}
L = K \oplus K \oplus \cdots \oplus K
\end{equation}
be the direct sum of  $k$ copies of $K$.  Then $L$ is a positive
definite even lattice with symmetric bilinear form given by
\begin{equation}
\langle (\alpha_1, \alpha_2,\dots,\alpha_k), (\beta_1, \beta_2, \dots ,
\beta_k ) \rangle = \sum_{j = 1}^k \langle \alpha_j, \beta_j \rangle,
\end{equation}
for $\alpha_j, \beta_j \in K$, $j = 1,\dots, k$.  The vertex operator algebra
associated to the lattice $L$ satisfies $V_L = V_K \otimes V_K \otimes
\cdots \otimes V_K = V_K^{\otimes k}$, where $V_K^{\otimes k}$ denotes
the $k$-fold tensor product of $V_K$ with itself.

Let $\nu \in \mathrm{Aut} \; L$ be given by
\begin{eqnarray}\label{define-nu}
\nu : K \oplus K \oplus \cdots \oplus K &\longrightarrow& K \oplus K \oplus
\cdots \oplus K \\
(\alpha_1, \alpha_2, \dots , \alpha_k) & \mapsto & (\alpha_2, \alpha_3, \dots,
\alpha_k, \alpha_1) .\nonumber
\end{eqnarray}
Then $\nu$ is an isometry of $L$, i.e., $\langle \nu\alpha,
\nu\beta \rangle = \langle \alpha,\beta \rangle$ for all $\alpha,
\beta \in L$.   As noted in Section \ref{lattice-section},  $\nu$ lifts
canonically to an automorphism $\hat{\nu}$ of $V_L = V_K^{\otimes k}$, and
in this setting, i.e., with $\nu$ given by (\ref{define-nu}), this automorphism
is
given by
\begin{eqnarray}\label{extend-nu}
\hat{\nu}: V_K \otimes V_K \otimes \cdots \otimes V_K & \longrightarrow & V_K
\otimes V_K \otimes \cdots \otimes V_K \\
v_1\otimes v_2 \otimes \dots \otimes v_k & \mapsto & v_2 \otimes v_3 \otimes
\dots \otimes v_k \otimes v_1 .\nonumber
\end{eqnarray}
That is, the automorphism is a ``permutation'' of $V_K^{\otimes k}$.
Thus it is appropriate to consider both the construction of
$\hat{\nu}$-twisted modules for the vertex operator algebra $V_L$ as
developed in \cite{L1} and \cite{FLM2} and the construction of
$\hat{\nu}$-twisted $V_L$-modules as developed in \cite{BDM}.

\begin{rema}{\em
In \cite{BDM}, the construction of $g$-twisted $V_L$-modules for $g$
any permutation on $k$ letters first relies on the construction of
$\hat{\nu}$-twisted $V_L$-modules for $\nu = (1 \; 2 \cdots k)$.  Thus we
first restrict ourselves to this particular permutation.  At the end of Section
\ref{explicit-section}, we discuss generalizations to arbitrary permutations. }
\end{rema}

\subsection{The ``space-time" construction of $\hat{\nu}$-twisted
$V_L$-modules}\label{spacetime-section}

Following the construction of twisted modules for a lattice and
isometry as developed in \cite{L1} and \cite{FLM2} (and see also
\cite{DL}) specialized to the setting introduced above, observe that
if $k$ is even, then for $\alpha = (\alpha_1, \dots, \alpha_k) \in L$,
we have
\begin{eqnarray}
\langle \nu^{k/2} \alpha, \alpha \rangle &=& \langle ( \alpha_{\frac{k}{2} + 1},
\alpha_{\frac{k}{2} + 2}, \dots, \alpha_k, \alpha_1, \alpha_2, \dots,
\alpha_{\frac{k}{2}}), (\alpha_1, \dots, \alpha_k) \rangle \\
&=& \langle \alpha_{\frac{k}{2}+1}, \alpha_1 \rangle + \langle
\alpha_{\frac{k}{2} + 2}, \alpha_2 \rangle + \cdots + \langle \alpha_k,
\alpha_{\frac{k}{2}} \rangle \nonumber \\
& & \quad + \; \langle \alpha_1, \alpha_{\frac{k}{2} + 1} \rangle + \langle
\alpha_2, \alpha_{\frac{k}{2} + 2} \rangle + \cdots + \langle
\alpha_{\frac{k}{2}}, \alpha_k \rangle \nonumber \\
&=& 2 \sum_{j = 1}^{k/2} \langle \alpha_j, \alpha_{\frac{k}{2} + j} \rangle .
\nonumber
\end{eqnarray}
Thus $\langle \nu^{k/2} \alpha, \alpha \rangle \in 2 \mathbb{Z}$ for $\alpha \in
L$, verifying Equation (\ref{even-equation}) in this setting.  This implies that
\begin{equation}\label{evenness}
\left< \sum_{j = 0}^{k-1} \nu^j \alpha, \alpha \right> \in 2 \mathbb{Z} ,
\end{equation}
verifying Equation (\ref{even-condition}) in this setting.  Thus the
commutator map $C$ given by (\ref{commutator-definition}) satisfies
$C(\alpha, \alpha) = 1$ for $\alpha \in L$.

Recalling our fixed primitive $k$-th root of unity $\eta$ {}from Section
\ref{lattice-section}, for $n \in \mathbb{Z}$ set
\begin{equation}\label{hngrading}
\h_{(n)} = \{ h \in \h \; | \; \nu h = \eta^{n} h \} \subset \h,
\end{equation}
so that $\h = \coprod_{n \in \mathbb{Z}/k\mathbb{Z}} \h_{(n)}$, where
we identify $\h_{(n \; \mathrm{mod} \; k)}$ with $\h_{(n)}$ for $n \in
\mathbb{Z}$. Then in general,
\begin{equation}
\h_{(n)} = \{ h + \eta^{-n} \nu h + \eta^{-2n} \nu^2 h + \cdots +
\eta^{-(k-1)n} \nu^{k-1} h \; | \; h \in \h \}
\end{equation}
and thus in the present setting
\begin{equation}
\h_{(n)} = \mathrm{span}_\C \{ (\alpha, \eta^{n} \alpha, \eta^{2n} \alpha,
\dots, \eta^{(k-1)n} \alpha) \; | \; \alpha \in K \},
\end{equation}
since
\begin{multline*}
(\alpha, \eta^{n} \alpha, \eta^{2n} \alpha, \dots, \eta^{(k-1)n} \alpha) =
(\alpha,0,\dots,0) + \eta^{n} (0,\alpha, 0,\dots, 0) \\
+ \eta^{2n}(0,0,\alpha,0,\dots,0) + \cdots + \eta^{(k-1)n} (0,\dots,0,\alpha) .
\end{multline*}

For $n \in \mathbb{Z}/k\mathbb{Z}$, denote by
\begin{equation}\label{Pn}
P_n : \h \longrightarrow \h_{(n)}
\end{equation}
the projection onto $\h_{(n)}$, and for $h \in \h$ and $n \in
\mathbb{Z}$, set $h_{(n)} = P_{(n \; \mathrm{mod}\; k)} h$.  In
general, we have that for $h \in \h$ and $n \in \Z$,
\begin{equation}\label{h_n-formula}
h_{(n)} = \frac{1}{k} \sum_{j=0}^{k-1} \eta^{-nj} \nu^j h,
\end{equation}
so that in the present setting, for $(\alpha_1, \alpha_2, \dots,
\alpha_k) \in L$,
\begin{equation}\label{alpha_n}
(\alpha_1, \alpha_2, \dots, \alpha_k) _{(n)} = \frac{1}{k} \Bigl( \sum_{j = 1}^k
\eta^{n(1-j)} \alpha_j, \sum_{j = 1}^k \eta^{n(2-j)} \alpha_j, \dots, \sum_{j =
1}^k \eta^{n(k-j)} \alpha_j \Bigr) .
\end{equation}

Viewing $\h$ as an abelian Lie algebra, consider the $\nu$-twisted
affine Lie algebra
\begin{equation}
\hat{\h}[\nu] = \coprod_{n\in\frac{1}{k} \Z} \h_{(kn)}\otimes
t^{n}\oplus \C {\bf c}
\end{equation}
with brackets determined by
\begin{equation}
[\alpha \otimes t^m, \beta \otimes t^n]=\langle \alpha , \beta \rangle
m\delta_{m+n,0}{\bf c}
\end{equation}
for $\alpha \in{\h}_{(km)}$, $\beta \in{\h}_{(kn)}$, and
$m,n\in\frac{1}{k} \Z$, and
\begin{equation}
[{\bf c},\hat{\h}[\nu]]=0.
\end{equation}
Define the {\em weight gradation} on $\hat{\h}[\nu]$ by
\begin{equation}\label{grade-h}
{\rm wt}\,(\alpha \otimes t^{n})=-n, \ \ \ {\rm wt}\,{\bf c}=0
\end{equation}
for $n\in \frac{1}{k} \Z$, $\alpha \in {\h}_{(kn)}$.  Set
\begin{equation}
\hat{\h}[\nu]^+=\coprod_{n>0}{\h}_{(kn)}\otimes t^{n},\ \ \quad
\hat{\h}[\nu]^-=\coprod_{n<0}{\h}_{(kn)}\otimes t^{n}.
\end{equation}
Now the subalgebra
\begin{equation}
\hat{\h}[\nu]_{\frac{1}{k} \Z}=\hat{\h}[\nu]^+\oplus
\hat{\h}[\nu]^-\oplus \C {\bf c}
\end{equation}
of $\hat{\h}[\nu]$ is a Heisenberg algebra. Form the induced
$\hat{\h}[\nu]$-module
\begin{equation}
S[\nu]=U(\hat{\h}[\nu])\otimes_{U(\coprod_{n \ge 0}{\h}_{(kn)}\otimes
t^{n}\oplus \C {\bf c})} \C \simeq S(\hat{\h}[\nu]^-)\ \ \ {\rm (linearly)},
\end{equation}
where $\coprod_{n\ge 0} \h_{(kn)} \otimes t^{n}$ acts trivially on
$\C$ and ${\bf c}$ acts as 1.  Then $S[\nu]$ is irreducible under
$\hat{\h}[\nu]_{\frac{1}{k} \Z}.$

Following \cite{DL}, Section 6, we give the module $S[\nu]$ the
natural $\Q$-grading (by weights) compatible with the action of
$\hat{\h}[\nu]$ and such that
\begin{eqnarray}\label{grade-vacuum}
{\rm wt}\,1 &=& \frac{1}{4k^2} \sum_{j = 1}^{k-1} j (k-j) {\rm
dim}\,({\h}_{(j)}) \\
&=& \frac{d}{4k^2} \sum_{j = 1}^{k-1} j (k-j), \nonumber
\end{eqnarray}
where
\begin{equation}
d = {\rm rank} \, K.
\end{equation}
Now
\begin{equation}\label{weight-identity}
\sum_{j = 1}^{k-1} j (k-j) = \frac{k(k^2-1)}{6}
\end{equation}
since, for example, $\sum j (k-j) = k \sum j - \sum j^2$, and thus
(\ref{grade-vacuum}) simplifies to
\begin{equation}\label{grade-vacuum2}
{\rm wt}\,1 =  \frac{(k^2-1)d}{24k}.
\end{equation}
Later we will justify (\ref{grade-vacuum2}) by determining the action
of the operator $L^{\hat{\nu}}(0)$ obtained {}from the general twisted vertex
operators introduced in \cite{FLM2}.

Following Sections 5 and 6 of \cite{L1} implemented in this special case,
we have that the automorphisms of $\Lnu$ covering the identity
automorphism of $L$ are precisely the maps $\rho^* : a \rightarrow a \rho
(\bar{a})$ for a homomorphism $\rho : L \rightarrow \langle \eta_0 \rangle$.
We have that
\begin{equation}\label{diagonallattice}
L \cap \h_{(0)} = \{ (\alpha, \alpha,\dots, \alpha) \; | \; \alpha \in
K \},
\end{equation}
the ``diagonal'' lattice, and there is a homomorphism $\rho_0 : L \cap
\h_{(0)} \rightarrow \langle \eta_0 \rangle$ such that $\hat{\nu} a =
a \rho_0(\bar{a})$ if $\nu \bar{a} = \bar{a}$.  Now $\rho_0$ can be
extended to a homomorphism $\rho : L \rightarrow \langle \eta_0
\rangle$ since the map $1 - P_0$ induces an isomorphism {}from
$L/L\cap \h_{(0)}$ to the free abelian group $(1 - P_0)L$.
Multiplying $\hat{\nu}$ by the inverse of $\rho^*$ gives us an
automorphism $\hat{\nu}$ of $\Lnu$ satisfying (\ref{nu-on-hatL}) and
\begin{equation}\label{propertyofnuhat}
\hat{\nu} a = a \quad \mathrm{if} \quad \nu \bar{a} = \bar{a},
\end{equation}
as in (\ref{nuhata=a}).

Let
\begin{equation}
N = (1-P_0)\h \cap L = \{\alpha \in L \; | \; \langle \alpha,
\h_{(0)} \rangle= 0 \}.
\end{equation}
Then
\begin{equation}
N= \{ (\alpha_1, \alpha_2, \dots, \alpha_{k-1}, - \alpha_1 - \alpha_2 - \cdots -
\alpha_{k-1}) \; | \; \alpha_j \in K, \; j=1,\dots,k \} .
\end{equation}
Let
\begin{equation}
M = (1-\nu)L \subset N.
\end{equation}
Then
\begin{eqnarray*}
M &=& \{ (\alpha_1, \alpha_2, \dots, \alpha_k) - (\alpha_2, \alpha_3,\dots,
\alpha_k, \alpha_1) \; | \; \alpha_j \in K, \; j = 1,\dots, k \} \\
&=& \{ (\alpha_1 - \alpha_2, \alpha_2 - \alpha_3,\dots, \alpha_{k-1} - \alpha_k,
\alpha_k - \alpha_1) \; | \; \alpha_j \in K, \; j=1,\dots,k \} \\
&=& N.
\end{eqnarray*}

For $\alpha \in \h$, we have $\sum_{j=0}^{k-1} \nu^j \alpha \in
\h_{(0)}$ and thus for $\alpha, \beta \in N$, the commutator map $C$,
defined by (\ref{commutator-definition}), simplifies to
\begin{equation}\label{C-on-N1}
C_N ( \alpha, \beta) = \eta^{ \sum_{j=0}^{k-1} \langle j \nu^j \alpha, \beta
\rangle} .
\end{equation}
We further find that for $\alpha_j, \beta_j \in K$ ($j = 1,\dots,
k-1$), we have
\begin{eqnarray*}
\lefteqn{\sum_{j=0}^{k-1} \langle j \nu^j (\alpha_1, \alpha_2, \dots,
\alpha_{k-1}, - \alpha_1 - \alpha_2 - \cdots - \alpha_{k-1}), (\beta_1,
\beta_2, \dots, \beta_{k-1}, - \beta_1 - \beta_2 } \\
& & \quad - \cdots - \beta_{k-1}) \rangle \\
&=& \sum_{j=1}^{k-1} \langle j (\alpha_{j+1} , \alpha_{j+2}, \dots,
\alpha_{k-1}, - \alpha_1 - \alpha_2 - \cdots - \alpha_{k-1}, \alpha_1,
\alpha_2, \dots, \alpha_j), (\beta_1, \\
& & \quad \beta_2, \dots, \beta_{k-1}, - \beta_1 - \beta_2 - \cdots -
\beta_{k-1}) \rangle \\
&=& \sum_{j=1}^{k-1} j\Bigl( \langle \alpha_{j+1} , \beta_1 \rangle + \langle
\alpha_{j+2}, \beta_2 \rangle + \cdots + \langle \alpha_{k-1}, \beta_{k-j-1}
\rangle + \langle - \alpha_1 - \alpha_2 - \cdots \\
& & \quad - \; \alpha_{k-1}, \beta_{k-j} \rangle + \langle \alpha_1,
\beta_{k-j+1} \rangle + \langle \alpha_2,  \beta_{k-j+2} \rangle + \cdots +
\langle \alpha_{j-1} , \beta_{k-1} \rangle \\
& & \quad + \;  \langle \alpha_j, - \beta_1 - \beta_2 - \cdots - \beta_{k-1})
\rangle \Bigr) \\
&=& \sum_{j=1}^{k-1} j\Bigl( \langle \alpha_{j+1} , \beta_1 \rangle + \langle
\alpha_{j+2}, \beta_2 \rangle + \cdots + \langle \alpha_{k-1}, \beta_{k-j-1}
\rangle - \langle \alpha_1, \beta_{k-j} \rangle \\
& & \quad - \; \langle \alpha_2, \beta_{k-j} \rangle - \cdots - \langle
\alpha_{k-1}, \beta_{k-j} \rangle + \langle \alpha_1, \beta_{k-j+1} \rangle +
\langle \alpha_2, \beta_{k-j+2} \rangle + \cdots \\
& & \quad + \; \langle \alpha_{j-1} , \beta_{k-1} \rangle - \langle \alpha_j,
\beta_1 \rangle - \langle \alpha_j , \beta_2 \rangle - \cdots - \langle
\alpha_j, \beta_{k-1} \rangle \Bigr) \\
&=& \sum_{n=1}^{k-1} \left< \alpha_n, \sum_{j=1}^{n-1} (n-j) \beta_j - n
\sum_{j=1}^{k-1} \beta_j + \sum_{j = n+1}^{k-1} (k-j+n) \beta_j -
\sum_{j=1}^{k-1} (k-j) \beta_j \right> \\
&=& \sum_{n=1}^{k-1} \left< \alpha_n, \sum_{j=1}^{n-1} (n-j-n-(k-j)) \beta_j +
(-n-(k-n)) \beta_n \right. \hspace{.7in}\\
& & \left. \quad + \sum_{j = n+1}^{k-1} (-n+k-j+n-(k-j)) \beta_j \right> \\
&=& \sum_{n=1}^{k-1} \left< \alpha_n, - \sum_{j=1}^{n} k \beta_j \right>\\
&=& - k \sum_{n=1}^{k-1} \sum_{j=1}^n \langle \alpha_n, \beta_j \rangle .
\end{eqnarray*}
Thus, on $N$, the commutator map $C$ further simplifies {}from Equation
(\ref{C-on-N1}) to
\begin{equation}\label{C-radical}
C_N ( \alpha, \beta ) = 1 .
\end{equation}
Let
\begin{equation}
R = \{\alpha \in N \; | \; C_N(\alpha, N) = 1 \}
\end{equation}
denote the radical of $C_N$, so that {}from (\ref{C-radical}), we have
\begin{equation}
R = N = M.
\end{equation}

Continuing to follow \cite{L1}, we denote by $\hat{Q}$ the subgroup of
$\Lnu$ obtained by pulling back any subgroup $Q$ of $L$.  Then
\begin{equation}
\hat{N} = \hat{M} = \hat{R} \simeq N \times \langle \eta_0 \rangle,
\end{equation}
an abelian group.  Observe that $a\hat{\nu}a^{-1} \in \hat{M} =
\hat{N}$ for all $a \in \Lnu$.  By Proposition 6.1 of \cite{L1}, there
exists a unique homomorphism $\tau: \hat{M} = \hat{N} \rightarrow
\C^\times$ such that
\begin{equation}
\tau(\eta_0) = \eta_0 \ \ \mbox{and}\ \
\tau(a\hat{\nu}a^{-1}) = \eta^{-\sum_{j=0}^{k-1} \langle \nu^j \bar{a}, \bar{a}
\rangle/2}
=\eta^{-k \langle\bar{a}_{(0)},\bar{a}_{(0)}\rangle/2}
\end{equation}
for $a \in \Lnu$ (recall (\ref{evenness})).  Denote by $\C_\tau$ the
one-dimensional $\hat{N}$-module $\C$ with character $\tau$ and write
\begin{equation}
T = \C_\tau;
\end{equation}
this is the (unique up to equivalence) irreducible $\hat{N}$-module
given by Proposition 6.2 of \cite{L1}.

Form the induced $\Lnu$-module
\begin{equation}
U_T = \C[\Lnu] \otimes_{\C[\hat{N}]} T \simeq \C[L/N].
\end{equation}
Then $\Lnu$ and $\h_{(0)}$ act on $U_T$ as follows:
\begin{eqnarray}
a \cdot b \otimes r &=& ab \otimes r, \\
h \cdot b \otimes r &=& \langle h, \bar{b} \rangle b \otimes r
\end{eqnarray}
for $a,b \in \Lnu$, $r\in T=\C_\tau$, $h \in \h_{(0)}$.  As operators on $U_T$,
\begin{equation}\label{h-operates}
ha=a(\langle h,\bar a\rangle+h)
\end{equation}
for $a \in \Lnu$ and $h \in \h_{(0)}$.  Since the projection map $P_0$
(recall (\ref{Pn})) induces an isomorphism {}from $L/N$ onto $P_{0}L,$
we have
\begin{equation}
U_T = \C[P_{0}L],
\end{equation}
and since
\begin{equation}
P_{0}L = \frac{1}{k} \left( L \cap \h_{(0)} \right)
\end{equation}
(recall (\ref{diagonallattice})), we have
\begin{equation}
U_T \simeq \C\left[\frac{1}{k} \left(L \cap \h_{(0)}\right)\right].
\end{equation}

\begin{rema}\label{lattice-quotient-remark} {\em
Therefore we have that $\C[P_0L]$ (and thus $U_T$) is isomorphic to
$\C[K]$ by extension of the isomorphism
\begin{eqnarray*}
f : P_0L & \longrightarrow & K\\
\frac{1}{k} (\alpha, \alpha,\dots, \alpha) & \mapsto & \alpha,
\end{eqnarray*}
for $\alpha \in K$. }
\end{rema}

Note that we can write
\begin{equation}
U_T=\coprod_{\alpha \in P_{0}L} U_\alpha,
\end{equation}
where
\begin{equation}
U_\alpha =\{u \in U_T \; | \; h \cdot u = \langle h, \alpha \rangle u \ \ {\rm
for \ \ } h \in \h_{(0)}\},
\end{equation}
and
\begin{equation}
a \cdot U_\alpha \subset U_{\alpha +\bar a_{(0)}}
\end{equation}
for $a\in \Lnu$ and $\alpha \in P_{0}L.$

We define an End\,$U_T$-valued formal Laurent series $x^h$ for $h\in
\h_{(0)}$ as follows:
\begin{equation}
x^h\cdot u = x^{\langle h,\alpha \rangle}u \ \ \ {\rm for}\ \ \alpha \in P_{0}L
\ \ {\rm and} \ \ u \in U_\alpha.
\end{equation}
Then {}from (\ref{h-operates}),
\begin{equation}
x^ha=ax^{\langle h,\bar a\rangle+h}\ \ \ {\rm for}\ \ \ a\in \Lnu
\end{equation}
as operators on $U_T.$ Also, for $h \in h_{(0)},$ if $\langle h, L
\rangle\in \Z$, define the operator $\eta^h$ on $U_T$ by
\begin{equation}\label{oct1}
\eta^h \cdot u = \eta^{\langle h,\alpha\rangle} u
\end{equation}
for $u \in U_{\alpha}$ with $\alpha \in P_{0}L$.  Then for
$a\in \Lnu$, we have
\begin{equation}\label{right}
\hat\nu a =a \eta^{-\sum_{j=0}^{k-1} \nu^j \bar{a} -\sum_{j=0}^{k-1}
\langle \nu^j \bar{a}, \bar{a} \rangle/2} =a \eta^{-k\bar a_{(0)} - k
\langle\bar a_{(0)},\bar a_{(0)}\rangle/2}
\end{equation}
as operators on $U_T$.

Define a $\C$-gradation on $U_T$ by
\begin{equation}\label{grade-U}
{\rm wt}\, u = \frac{1}{2}\langle \alpha , \alpha \rangle \ \ \ {\rm
for}\ \ \alpha \in P_{0}L \ \ {\rm and} \ \ u \in U_\alpha.
\end{equation}

Form the space
\begin{eqnarray}
V^T_L &=& S[\nu]\otimes U_T \\ &=&
\left(U(\hat{\h}[\nu])\otimes_{U(\coprod_{n \ge 0}{\h}_{(kn)}\otimes
t^{n}\oplus \C {\bf c})} \C\right)
\otimes \left(\C[\Lnu]  \otimes_{\C[\hat{N}]} \C_\tau\right) \nonumber \\
&\simeq & S(\hat{\h}[\nu]^-) \otimes \C[P_{0}L], \nonumber
\end{eqnarray}
which is naturally graded (by weights), using the weight gradations of
$S[\nu]$ and $U_T.$

We let $\Lnu,$ $\hat{\h}[\nu]_{\frac{1}{k} \Z},$ ${\h}_{(0)}$ and
$x^h$ for $h\in{\h}_{(0)}$ act on $V_L^T$ by acting on either $S[\nu]$
or $U_T$, as described above.

For $\alpha \in \h$ and $n\in \frac{1}{k} \Z$, write $\alpha^T (n)$ or
$\alpha_{(kn)} (n)$ for the operator on $V_L^T$ associated with
$\alpha_{(kn)}\otimes t^n$, and set
\begin{equation}
\alpha^T (x)=\sum_{n\in\frac{1}{k} \Z} \alpha^T (n)x^{-n-1}
=\sum_{n\in\frac{1}{k} \Z} \alpha_{(kn)} (n) x^{-n-1}.
\end{equation}
Note that for $\alpha = (\alpha_1,\alpha_2, \dots, \alpha_k) \in \h$,
{}from (\ref{alpha_n}) we have
\begin{eqnarray*}
\alpha^T (x) &=& \sum_{n \in \frac{1}{k} \Z} \alpha^T (n) x^{-n-1} \\
&=& \sum_{n \in \frac{1}{k} \Z} \frac{1}{k} \Bigl( \sum_{j = 1}^k \eta^{kn(1-j)}
\alpha_j, \sum_{j = 1}^k \eta^{kn(2-j)} \alpha_j, \dots, \sum_{j = 1}^k
\eta^{kn(k-j)} \alpha_j \Bigr)(n) x^{-n-1} \\
&=& \frac{1}{k} \sum_{n \in \Z} \Bigl( \sum_{j = 1}^k \eta^{n(1-j)} \alpha_j,
\sum_{j = 1}^k \eta^{n(2-j)} \alpha_j, \dots, \sum_{j = 1}^k \eta^{n(k-j)}
\alpha_j \Bigr)\Bigl(\frac{n}{k} \Bigr) x^{-n/k-1}.
\end{eqnarray*}

Following \cite{L1} and \cite{FLM2}, for $\alpha \in L$, define
\begin{equation}
\sigma(\alpha) = \left\{ \begin{array}{ll}
\displaystyle{ \prod_{0< j < k/2} (1-\eta^{-j})^{\langle \nu^j \alpha, \alpha
\rangle} } 2^{\langle \nu^{k/2} \alpha, \alpha \rangle/2}& \mbox{if $k \in
2\mathbb{Z}$} \\
\\
\displaystyle{\prod_{0 < j < k/2  } (1-\eta^{-j})^{\langle \nu^j \alpha, \alpha
\rangle} } & \mbox{if $k \in 2\mathbb{Z} + 1$.}\\
\end{array}
\right.
\end{equation}
Then $\sigma(\nu\alpha)=\sigma(\alpha)$.  Using the normal-ordering
procedure described above, define the {\it $\hat{\nu}$-twisted vertex
operator} $Y^{\hat{\nu}}(a,x)$ for $a\in \hat{L}$ acting on $V_L^T$ as
follows:
\begin{equation}\label{L-operator}
Y^{\hat{\nu}}(a,x)= k^{-\langle \bar{a},\bar{a}\rangle /2} \sigma(\bar{a}) \
_\circ^\circ
e^{\int(\bar{a}^T (x)-\bar{a}(0)x^{-1})} a x^{\bar{a}_{(0)}+\langle
\bar{a}_{(0)}
,\bar{a}_{(0)}\rangle /2-\langle \bar{a},\bar{a}\rangle /2} \ _\circ^\circ .
\end{equation}
Note that on the right-hand side of (\ref{L-operator}), we view $a$ as an
element of $\Lnu$ using our set-theoretic identification between $\hat{L}$
and $\Lnu$ given by (\ref{identify-central-extensions}).
For $\alpha_1,\dots,\alpha_m \in{\h},$ $n_1,\dots,n_m \in \Z_+$ and
$v=\alpha_1(-n_1)\cdots \alpha_m(-n_m) \cdot \iota(a)\in V_L$, set
\begin{multline}
W(v,x)= \ _\circ^\circ
\biggl(\frac{1}{(n_1-1)!}\left(\frac{d}{dx}\right)^{n_1-1}\alpha_1^T(x)\biggr)\cdots
\\
\left(\frac{1}{(n_m-1)!}\biggl(\frac{d}{dx}\right)^{n_m-1}
\alpha_m^T(x)\biggr)Y^{\hat{\nu}}(a,x) \ _\circ^\circ ,
\end{multline}
where the right-hand side is an operator on $V^T_L$. Extend to all
$v\in V_L$ by linearity.

Define constants $c_{mnr} \in \C$ for $m, n \in \mathbb{N}$ and $r =
0,\dots, k-1$ by the formulas
\begin{eqnarray}\label{define-c's1}
\sum_{m,n\ge 0} c_{mn0} x^m y^n &=& -\frac{1}{2} \sum_{j = 1}^{k-1}
{\rm log} \left(\frac {(1+x)^{1/k} - \eta^{-j}
(1+y)^{1/k}}{1-\eta^{-j} }\right),\\
\sum_{m,n\ge 0} c_{mnr} x^m y^n &= & \frac{1}{2}{\rm log} \left( \frac
{(1+x)^{1/k} -\eta^{-r}
(1+y)^{1/k}}{1-\eta^{-r}}\right) \ \ \mbox{for}\ \ r \ne 0. \label{define-c's2}
\end{eqnarray}
(These are well-defined formal power series in $x$ and $y$.)  Let
$\{\beta_1,\dots, \beta_{\dim \h}\}$ be an orthonormal basis of $\h$,
and set
\begin{equation}\label{define-Delta_x}
\Delta_x = \sum_{m,n\ge 0} \sum_{r=0}^{k-1} \sum^{\dim \h}_{j=1}
c_{mnr} (\nu^{-r} \beta_j)(m) \beta_j(n) x^{-m-n} .
\end{equation}
Then $e^{\Delta_x}$ is well defined on $V_L$ since $c_{00r}=0$ for all
$r$, and for $v\in V_L,$ $e^{\Delta_x}v\in V_L[x^{-1}]$.  Note that
$\Delta_x$ is independent of the choice of orthonormal basis.  In our
special case, recall that $d = {\rm rank} \, K$, so that ${\rm dim} \,
\h = kd$.

For $v\in V_L,$ the $\hat{\nu}$-$twisted$ $vertex$ $operator$
$Y^{\hat{\nu}}(v,x)$ is
defined by:
\begin{equation}\label{Ynuhat}
Y^{\hat{\nu}}(v,x)=W(e^{\Delta_x}v,x).
\end{equation}
Then this yields a well-defined linear map
\begin{eqnarray}
V_L &\longrightarrow&(\mbox{End}\,V^T_L)[[x^{1/k},x^{-1/k}]] \\ \ v
&\mapsto& Y^{\hat{\nu}}(v,x)= \sum_{n \in \frac{1}{k}\Z}v^{\hat{\nu}}_nx^{-n-1}
\nonumber
\end{eqnarray}
where $v^{\hat{\nu}}_n\in {\rm End}\,V^T_L$.  Recall {}from
(\ref{nuhat^k=1}) that $\hat{\nu}$ has period (and hence order) $k$ on
$\hat{L}$, and thus on the vertex operator algebra $V_L$ as well.

It has been established in \cite{L1}, \cite{FLM2}, \cite{L2} and \cite{DL}
that $(V_L^T, Y^{\hat{\nu}})$ is an irreducible $\hat{\nu}$-twisted
$V_L$-module (recall Section \ref{generalities-section} for the
definition).

Now, following \cite{DL} (but filling in some details), we will
justify the weight gradation of $V_L^T$  given by (\ref{grade-h}),
(\ref{grade-vacuum}) (or equivalently, (\ref{grade-vacuum2})), and
(\ref{grade-U}) by showing that this grading is given by the
eigenvalues of the operator $L^{\hat{\nu}}(0)$ (recall (\ref{define-L^g}) and
(\ref{L^g-grading})).  This will then allow us to calculate the graded
dimension of $V_L^T$, which is the main piece of data we will use to
establish an isomorphism between the space-time $\hat{\nu}$-twisted
$V_L$-module construction just established above and the worldsheet
$\hat{\nu}$-twisted $V_L$-module construction given in Section
\ref{bdm-section}.

We first note that for $\alpha, \beta \in \mathfrak{h}$ and $s,t =
0,\dots, k-1$, we have, by (\ref{hngrading}),
\begin{equation*}
\langle \alpha_{(s)}, \beta_{(t)} \rangle
= \langle \nu \alpha_{(s)}, \nu \beta_{(t)} \rangle
= \eta^{s+t} \langle \alpha_{(s)}, \beta_{(t)} \rangle,
\end{equation*}
so that
\begin{equation}\label{r=-s}
\langle \alpha_{(s)}, \beta_{(t)} \rangle = 0
\ \ \mbox{unless}\ \ s+t \equiv 0 \ \ \mbox{mod}\,k.
\end{equation}

{}From  (\ref{define-Delta_x}) and (\ref{r=-s}), we have
\begin{eqnarray*}
\lefteqn{\Delta_x \cdot \alpha (-1) \beta (-1) {\bf 1}}\\
&=&  \sum_{m,n\ge 0}  \sum_{r=0}^{k-1} \sum^{\dim \h}_{j=1} c_{mnr} (\nu^{-r}
\beta_j)(m) \beta_j(n) x^{-m-n}  \cdot \alpha(-1) \beta (-1) {\bf 1} \\
&=&  \sum_{r=0}^{k-1} \sum^{\dim \h}_{j=1} c_{11r} (\nu^{-r} \beta_j)(1)
\beta_j(1)
x^{-2}  \cdot \alpha (-1) \beta (-1){\bf 1} \\
&=&  \sum_{r=0}^{k-1} \sum^{\dim \h}_{j=1} c_{11r} ( \langle \beta_j ,   \alpha
\rangle  \langle \beta_j , \nu^r \beta \rangle  {\bf 1} + \langle \beta_j,
\beta \rangle   \langle  \beta_j,  \nu^r \alpha \rangle {\bf 1} ) x^{-2}   \\
&=&  \sum_{r=0}^{k-1} c_{11r} ( \langle  \alpha, \nu^r \beta \rangle  {\bf 1} +
\langle \beta,  \nu^r \alpha \rangle {\bf 1} ) x^{-2}   \\
&=& \sum_{r=0}^{k-1} c_{11r} \sum_{s,t = 0}^{k-1} ( \langle  \alpha_{(s)}, \nu^r
\beta_{(t)} \rangle  {\bf 1} + \langle \beta_{(t)},  \nu^r \alpha_{(s)} \rangle
{\bf 1} ) x^{-2} \hspace{.7in} \\
&=& \sum_{r=0}^{k-1} c_{11r} \sum_{s,t = 0}^{k-1} (\eta^{rt} + \eta^{r s})
\langle  \alpha_{(s)}, \beta_{(t)} \rangle  {\bf 1}  x^{-2}   \\
&=& \sum_{r=0}^{k-1} c_{11r} \sum_{s= 0}^{k-1} (\eta^{-rs} + \eta^{rs})
\langle  \alpha_{(s)}, \beta_{(-s)} \rangle  {\bf 1} x^{-2}   .
\end{eqnarray*}
Thus
\begin{multline}\label{6.20}
e^{\Delta_x}\alpha(-1)\beta(-1) {\bf 1}=\alpha(-1)\beta(-1) {\bf 1} \\
+ \left(2c_{110} \langle\alpha,\beta\rangle  +\sum_{r=1}^{k-1} \sum_{s =
0}^{k-1} c_{11r}(\eta^{rs}+ \eta^{-rs})\langle  \alpha_{(s)}, \beta_{(-s)}
\rangle \right) {\bf 1} x^{-2} .
\end{multline}

For $s = 0,\dots, k-1$, let $\{\beta^{(s)}_1,\dots,\beta^{(s)}_{\dim
\h_{(s)}}\}$ be a basis of $\h_{(s)}$, and let
\begin{equation}
\{(\beta^{(s)}_{j_s})^* \ | \ j_s= 1,\dots, \dim  \h_{(s)}, \ s = 0, \dots, k-1
\}
\end{equation}
be a dual basis for $\h$ with respect to $\langle \cdot , \cdot \rangle$.  Then
$\langle (\beta^{(s)}_{j_s})_{(s)}, (\beta^{(s)}_{j_s})^*_{(-s)} \rangle =
\langle \beta^{(s)}_{j_s} , (\beta^{(s)}_{j_s})^* \rangle = 1$.  Recalling
(\ref{define-omega}), we have (cf. \cite{FLM3})
\begin{equation}
\omega = \frac{1}{2} \sum_{s= 0}^{k-1} \sum_{j_s = 1}^{\dim \h_{(s)}}
\beta_{j_s}^{(s)}(-1) (\beta_{j_s}^{(s)})^* (-1) \mathbf{1}
\end{equation}
and
\begin{eqnarray}\label{dims-not-equal}
\qquad e^{\Delta_x} \omega &=& \omega +  \frac{1}{2} \sum_{s = 0}^{k-1}
\sum_{j_s= 1}^{\dim \h_{(s)}}  \biggl(2c_{110} \langle \beta^{(s)}_{j_s},
(\beta^{(s)}_{j_{s}})^* \rangle \\
& & \quad +\sum_{r=1}^{k-1} c_{11r}(\eta^{rs}+ \eta^{-rs})\langle
\beta^{(s)}_{j_s}, (\beta^{(s)}_{j_s})^* \rangle \biggr) {\bf 1} x^{-2}
\nonumber\\
&=& \omega + \frac{1}{2}  \sum_{s = 0}^{k-1}   \biggl(2c_{110} \dim \h_{(s)}
+\sum_{r=1}^{k-1} c_{11r}(\eta^{rs}+ \eta^{-rs}) \dim \h_{(s)} \biggr) {\bf 1}
x^{-2} \nonumber\\
&=& \omega +    \biggl(c_{110} \dim \h +  \frac{1}{2} \sum_{s = 0}^{k-1}
\sum_{r=1}^{k-1} c_{11r}(\eta^{rs}+ \eta^{-rs}) \dim \h_{(s)} \biggr) {\bf 1}
x^{-2} .\nonumber
\end{eqnarray}
If $\dim \h_{(s)} = \dim \h_{(t)}$ for all $s,t = 0,\dots, k-1$ as in our
specialized setting, then (\ref{dims-not-equal}) further simplifies to
\begin{eqnarray}\label{dims-equal}
e^{\Delta_x} \omega &=& \omega +    \biggl(c_{110} \dim \h +  \frac{\dim \h}{2k}
\sum_{s = 0}^{k-1} \sum_{r=1}^{k-1} c_{11r}(\eta^{rs}+ \eta^{-rs})  \biggr) {\bf
1} x^{-2} \\
&=& \omega +  c_{110} \dim \h {\bf 1} x^{-2}. \nonumber
\end{eqnarray}

The number $c_{110}$ is defined by (\ref{define-c's1}) and can be expressed as
\begin{eqnarray}
c_{110} &=& \left. -\frac{\partial}{\partial x} \frac{\partial}{\partial y}
\frac{1}{2}\sum_{j = 1}^{k-1} {\rm log}\left(\frac {(1+x)^{1/k} - \eta^{-j}
(1+y)^{1/k}}{1-\eta^{-j} }\right) \right|_{x=y=0} \label{calculate-c_110} \\
&=&  -\frac{1}{2k^2}\sum_{j=1}^{k-1}\frac{\eta^{-j}} {(1-\eta^{-j})^2} \nonumber
\end{eqnarray}

The next lemma follows {}from Equations (6.21) and (6.22) in
\cite{DL}. Since the proof of this fact was not included in \cite{DL},
we supply it here for completeness.
\begin{lem}\label{mystery-lemma}
For any $m \in \mathbb{Z}_+$ and $\eta_m$ a primitive $m$-th root of
unity
\begin{equation} \label{mystery-equation}
\sum_{j = 1}^{m-1} \frac{\eta_m^{-j} }{(1-\eta_m^{-j})^2 } = - \frac{m^2 -1}
{12}
\end{equation}
\end{lem}
\begin{proof}
By direct expansion, we observe that
\begin{equation}\label{delta}
\frac{1}{m} \sum_{j \in \mathbb{Z}/m\mathbb{Z}} \delta(\eta_m^{-j} x) =
\delta(x^m) .
\end{equation}
Considering only the nonnegative powers of $x$ in (\ref{delta}), we have the
equality
\begin{equation}\label{before-D}
\frac{1}{m} \sum_{j \in \mathbb{Z}/m\mathbb{Z}} \frac{1}{1-\eta_m^{-j} x} =
\frac{1}{1-x^m}
\end{equation}
of formal rational functions, and applying $x \frac{d}{dx}$ to both sides of
(\ref{before-D}) gives
\begin{equation}
\frac{1}{m} \sum_{j \in \mathbb{Z}/m\mathbb{Z}}
\frac{\eta_m^{-j}x}{(1-\eta_m^{-j}
x)^2} =  \frac{mx^m}{(1-x^m)^2} .
\end{equation}
Therefore
\begin{equation}\label{rational-functions}
\sum_{j =1}^{m-1} \frac{\eta_m^{-j}x}{(1-\eta_m^{-j} x)^2} =
\frac{m^2x^m}{(1-x^m)^2} -   \frac{x}{(1-x)^2} .
\end{equation}
The left-hand side of (\ref{mystery-equation}) is obtained by setting
$x=1$ in the left-hand side of  (\ref{rational-functions}). To prove
(\ref{mystery-equation}), we take the limit as $x$ approaches 1 of the
right-hand side of (\ref{rational-functions}). An efficient method for
computing this limit is to replace $x$ by $x+1$ in the right-hand side
of (\ref{rational-functions}), and then take the limit as $x$
approaches 0.  Replacing $x$ by $x+1$ on the right-hand side of
(\ref{rational-functions}) and then dividing by $x+1$ (which
approaches 1 as $x$ approaches 0) gives
\begin{eqnarray*}
\lefteqn{ \frac{m^2(x+1)^{m-1}}{(1-(x+1)^m)^2} - \frac{1}{x^2} }\\
&=& \frac{m^2(x+1)^{m-1} - \left(\frac{1 - (x+1)^m}{x} \right)^2}{(1-(x+1)^m)^2}
\\
&=& \frac{ \sum_{n\in \mathbb{N}} m^2 \binom{m-1}{n} x^n - \left(\sum_{n \in
\mathbb{Z}_+} \binom{m}{n} x^{n-1} \right)^2}{\left( \sum_{n \in \mathbb{Z}_+}
\binom{m}{n} x^n \right)^2}  \\
&=& \frac{ m^2 + m^2(m-1)x +  m^2 \binom{m-1}{2} x^2 + O(x^3) - \left( m +
\binom{m}{2} x + \binom{m}{3} x^2 + O(x^3) \right)^2}{ \left(mx + O(x^2)
\right)^{2}} \\
&=& \frac{ m^2\binom{m-1}{2} x^2 - \binom{m}{2}^2 x^2 - 2 m \binom{m}{3}
x^2 + O(x^3) }{m^2x^2 + O(x^3)} \\
&=& \frac{1}{2} (m-1)(m-2) - \frac{1}{4} (m-1)^2 - \frac{1}{3} (m-1)(m-2) +
O(x) .\\
\end{eqnarray*}
Thus the limit as $x$ approaches 0 is
\begin{equation*}
(m-1) \left( \frac{6(m-2) - 3 (m-1) - 4(m-2)}{12} \right) = - \frac{m^2 -
1}{12},
\end{equation*}
proving (\ref{mystery-equation}).
\end{proof}

Thus we have that
\begin{equation}\label{c_110}
c_{110}  =  \frac{k^2 - 1}{24k^2}
\end{equation}
and
\begin{equation}
e^{\Delta_x} \omega = \omega +  \frac{k^2 - 1}{24k^2} \dim \h {\bf 1} x^{-2}.
\end{equation}
Note that of course this is independent of the choice of basis for $\h$.  Thus
for any orthonormal basis for $\h$, $\{\beta_1, \dots, \beta_{\dim \h} \}$, and
recalling  (\ref{define-omega}) and (\ref{Ynuhat}), we have
\begin{eqnarray}
Y^{\hat{\nu}}(\omega, x) &=& \frac{1}{2}\sum_{j=1}^{\dim \h} \ _\circ^\circ
\beta^T_j(x)\beta^T_j(x)_\circ^\circ + \frac{k^2 - 1}{24k^2} \dim \h  x^{-2}.
\\
&=& \sum_{n \in \mathbb{Z}} L^{\hat{\nu}}(n) x^{-n-2} .\nonumber
\end{eqnarray}
By the theorem quoted above that $V_L^T$ is a $\hat{\nu}$-twisted
$V_L$-module, the operators $L^{\hat{\nu}}(n)$ satisfy the Virasoro
algebra relations (\ref{Viralgrelations}).  As we now show, the
grading on $V_L^T$ described above is given by
$L^{\hat{\nu}}(0)$-eigenvalues.

Since in our case $\dim \h = kd$, we have
\begin{equation}\label{L-twisted}
L^{\hat{\nu}}(0) =  \frac{1}{2} \sum_{j = 1}^{kd} \sum_{n \in \frac{1}{k}
\mathbb{Z}}
\beta_j^T(-|n|) \beta_j^T(|n|) + \frac{(k^2 - 1)d}{24k}.
\end{equation}
Thus
\begin{equation}
L^{\hat{\nu}}(0)1 =  \frac{(k^2 - 1)d}{24k},
\end{equation}
as in (\ref{grade-vacuum2}).  Similarly, for $u=1\otimes u\in
V_L^T$ with $u\in U_{\alpha}\subset U_T$  ($\alpha\in  P_0L$), we have
\begin{eqnarray}
L^{\hat{\nu}}(0) u &=&  \frac{1}{2} \sum_{j = 1}^{kd} (\beta_j)_{(0)}
(\beta_j)_{(0)}
u  + \frac{(k^2 - 1)d}{24k} u  \\
&=&  \frac{1}{2} \sum_{j = 1}^{kd} \langle (\beta_j)_{(0)}, \alpha \rangle
\langle (\beta_j)_{(0)}, \alpha \rangle u  + \frac{(k^2 - 1)d}{24k} u \nonumber
\\
&=&  \left( \frac{1}{2}\langle\alpha,\alpha \rangle + \frac{(k^2 -
1)d}{24k}\right) u, \nonumber
\end{eqnarray}
and for $m\in \frac{1}{k} \mathbb{Z}$ and $\alpha \in \mathfrak{h}_{(km)}$
\begin{eqnarray}
[L^{\hat{\nu}}(0), \alpha^T(m) ] &=& - \sum_{j = 1}^{kd}   \beta_j^T (m) \langle
(\beta_j)_{(-km)} , \alpha_{(km)} \rangle m  \\
&=& - m \sum_{j = 1}^{kd}   \beta_j^T (m) \langle
\beta_j , \alpha_{(km)} \rangle \nonumber  \\
&=& - m \alpha^T(m) . \nonumber
\end{eqnarray}
Thus $L^{\hat{\nu}}(0) v = ({\rm wt} \; v + \frac{(k^2 - 1)d}{24k})v$ for $v \in
V_L^T$, using the
weight gradation defined by (\ref{grade-h}) and (\ref{grade-U}) and
incorporating the grading
shift given by (\ref{grade-vacuum2}).

Using this, we find that the graded dimension of the $\hat{\nu}$-twisted
$V_L$-module $V_L^T$ is
\begin{eqnarray}\label{character-twisted}
{\rm dim}_* V_L^T &=& {\rm tr}_{V_L^T} \; q^{L^{\hat{\nu}}(0) - kd/24} \nonumber
\\
&=& q^{ (k^2 - 1)d/ 24k - kd/24} \left(\sum_{\beta \in P_0L}
q^{ \langle \beta, \beta \rangle/2}
\right) \left( \prod_{n \in \Z_+} (1 - q^{n/k} )^{-d}\right) \nonumber\\
&=& q^{-d/24k} \left(\sum_{\alpha \in K}
q^{ \langle \frac{1}{k}(\alpha,\alpha,\dots, \alpha) ,
\frac{1}{k}(\alpha,\alpha,\dots, \alpha) \rangle/2}
\right) \left( \prod_{n \in \Z_+} (1 - q^{n/k} )^{-d}\right) \nonumber \\
&=& q^{ -d/24k} \left(\sum_{\alpha \in K} q^{ \langle \alpha , \alpha \rangle /
2k} \right) \left( \prod_{n \in \Z_+} (1 - q^{n/k})^{-d} \right) .
\end{eqnarray}

\section{The ``worldsheet" construction and classification of
$\hat{\nu}$-twisted $V_L$-modules}\label{bdm-section}

Following \cite{BDM} we give the construction of $\hat{\nu}$-twisted
$V_L$-modules for the case when $V_L = V_K^{\otimes k}$ for $K$ a
positive definite even lattice and for $\hat{\nu}$ given by
(\ref{extend-nu}).

Define $\mathcal{E}_f(x^{1/k})\in ({\rm End} \; V_K) [[x^{1/k}, x^{-1/k}]]$ by
\begin{equation}
\mathcal{E}_f (x^{1/k}) = \exp \Biggl( \sum_{j \in \Z_+} a_j x^{- j/k} L(j)
\Biggr) k^{-L(0)} x^{\left( 1/k -1 \right) L(0)}
\end{equation}
where the $L(j) \in {\rm End} \; V_K$, for $j \in \mathbb{N}$, are the
elements given by the vertex operator algebra structure on $V_K$ and
where the $a_j \in \C$, for $j \in \Z_+$, are given uniquely by
\begin{equation}
\exp \Biggl( - \sum_{j \in \Z_+} a_j  x^{j + 1} \frac{\partial}{\partial x}
\Biggr) \cdot x = \frac{1}{k} (1 + x)^k - \frac{1}{k} .
\end{equation}
For example, $a_1= (1-k)/2$ and $a_2= (k^2-1)/12.$

\begin{rema}\label{worldsheet-remark}{\em
We use the symbol $f$ in the operator  $\mathcal{E}_f(x^{1/k})$ and the term
``worldsheet" to describe the twisted construction we will recall
{}from \cite{BDM} for the following reason:
Let $w$, $y$ and $z$ be formal variables and consider the (formal)
function
\begin{equation}
f(y) = \exp \Biggl( - \sum_{j \in \Z_+} a_j z^{- j/k} y^{j+1}
\frac{\partial}{\partial y} \Biggr) k^{y \frac{\partial}{\partial y}}
z^{\left( 1-1/k \right) y \frac{\partial}{\partial y}} \cdot y .
\end{equation}
Then
\begin{eqnarray}
f(y) &=& z^{\left( -1/k  \right) y \frac{\partial}{\partial y}} \exp \Biggl( -
\sum_{j \in \Z_+} a_j y^{j+1} \frac{\partial}{\partial y} \Biggr) k^{y
\frac{\partial}{\partial y}}  z^{y \frac{\partial}{\partial y}}\cdot y \\
&=&  kz \Bigl( \frac{1}{k} (1 + z^{-1/k} y )^k - \frac{1}{k} \Bigr) \nonumber\\
&=& z(1 + z^{-1/k}y)^k - z \nonumber
\end{eqnarray}
and $f(y)$ has inverse $f^{-1}(y) = (y+z)^{1/k} - z^{1/k}$, which when
evaluated at $y = w-z$ gives $w^{1/k} - z^{1/k}$.  Now let $w$ and
$z$ be complex variables on the Riemann sphere and consider the
Riemann sphere with three punctures: at infinity, $z$ and zero.  Let
the local coordinate at $z$ be $w-z$.   Under the
correspondence between the geometry of propagating
strings and the algebra of vertex operators developed in \cite{H}, this
sphere with punctures corresponds to a certain ``worldsheet" -- a Riemann
surface swept out by propagating strings.   Choosing a branch
cut for the logarithm, we see that $f^{-1} (y)|_{y = w - z} = w^{1/k}
- z^{1/k}$ gives the local coordinate vanishing at $z^{1/k}$
on this three punctured sphere under the ``orbifolding" transformation
$w \mapsto w^{1/k}$.   The geometric interpretation of vertex
operator algebras developed in \cite{H} shows that
the operator corresponding to the change of variables
$f^{-1} (y)|_{y = w - z} = w^{1/k} - z^{1/k}$ in a vertex operator
algebra is $\mathcal{E}_f(z^{1/k})$. }
\end{rema}

\begin{rema} {\em
In \cite{BDM} the operator $\mathcal{E}_f(x^{1/k})$ is denoted
$\Delta_k (x)$.  We are changing notation to avoid confusion with the
notation $\Delta_x$ for the operator used to construct twisted modules
in Section \ref{specialize-section}. }
\end{rema}

For $v \in V_K$, define
\begin{equation}
v^1 = v \otimes \mathbf{1} \otimes \mathbf{1} \otimes \cdots \otimes
\mathbf{1} \in V_K^{\otimes k}
\end{equation}
and
\begin{equation}
v^{j+1} = \hat{\nu}^{-j} (v^1)
\end{equation}
for $j \in \Z$.  Thus $v^j$ is the element of $V_K^{\otimes k}$ that
has $v$ as the $(j \ \mathrm{mod} \ k)$-th tensor factor and
$\mathbf{1}$'s as the other tensor factors.

Let $(M, Y_K)$ be a $V_K$-module.  We will denote by $Y_{\hat{\nu}}$
the twisted operators on $M$ defined via the construction given in
\cite{BDM}, and they are defined as follows:
\begin{equation}\label{bdm-twisted1}
Y_{\hat{\nu}} ( u^1, x) \; = \; Y_K(\mathcal{E}_f (x^{1/k})u, x^{1/k})
\end{equation}
and
\begin{eqnarray}\label{bdm-twisted2}
Y_{\hat{\nu}} (u^{j+1}, x) &=& Y_{\hat{\nu}}(\hat{\nu}^{-j} (u^1), x) \\
&=& \lim_{x^{1/k} \rightarrow \eta^j x^{1/k}} Y_{\hat{\nu}}(u^1, x) \nonumber
\end{eqnarray}
for $u \in V_K$ and $\eta$ a fixed primitive $k$-th root of unity. Since
$V_K^{\otimes k}$ is generated by $u^j$ for $u \in V_K$ and $j = 1,\dots, k$,
the twisted vertex operators given in (\ref{bdm-twisted1}) and
(\ref{bdm-twisted2}) determine all the twisted vertex operators
$Y_{\hat{\nu}}(v, x)$ for $v \in V_L = V_K^{\otimes k}$.  In \cite{BDM}, it is
proved in particular that $(M, Y_{\hat{\nu}})$ is a $\hat{\nu}$-twisted
$V_K^{\otimes k}$-module and that $(M, Y_{\hat{\nu}})$ is irreducible if and
only if $(M, Y_K)$ is irreducible.

\begin{rema}{\em
In \cite{BDM}, the primitive $k$-th root of unity corresponding to $\eta$ is
fixed
to be $e^{2 \pi i/k}$.  However, the results of \cite{BDM} hold if $\eta$ is
chosen
to be any fixed primitive $k$-th root of unity.}
\end{rema}

On the other hand, letting $(M,Y_{\hat{\nu}})$ be a $\hat{\nu}$-twisted
$V_K^{\otimes k}$-module, we can define
\begin{equation}
Y^{\hat{\nu}}_K(u,x) = Y_{\hat{\nu}}((\mathcal{E}_f(x)^{-1} u)^1 , x^k),
\end{equation}
where
\begin{equation}
\mathcal{E}_f (x^{1/k})^{-1} = x^{\left( 1- 1/k \right) L(0)} k^{L(0)} \exp
\Biggl( -\sum_{j \in \Z_+} a_j x^{- j/k} L(j) \Biggr),
\end{equation}
and where we assume that if one replaces $x$ by a complex variable
$z$, then $z$ is restricted to complex values such that $(z^k)^{1/k} =
z$ for the standard branch cut of $\log$.  In \cite{BDM}, it is proved
in particular that $(M, Y^{\hat{\nu}}_K)$ is a $V_K$-module and that
$(M, Y^{\hat{\nu}}_K)$ is irreducible if and only if $(M,
Y_{\hat{\nu}})$ is irreducible.

Denote the category of $V_K$-modules by $\mathcal{ C}(V_K)$ and denote
the category of $\hat{\nu}$-twisted $V_L$-modules by $\mathcal{
C}^{\hat{\nu}}(V_L)$.  Define functors $F_{\hat{\nu}}$ and
$G_{\hat{\nu}}$ by
\begin{eqnarray*}
F_{\hat{\nu}} : \mathcal{ C}(V_K) &\longrightarrow& \mathcal{
C}^{\hat{\nu}}(V_L)\\
(M,Y_K) &\mapsto& (M,Y_{\hat{\nu}})\\
\end{eqnarray*}
and
\begin{eqnarray*}
G_{\hat{\nu}} : \mathcal{ C}^{\hat{\nu}} (V_L) &\longrightarrow& \mathcal{
C}(V_K)\\
(M,Y_{\hat{\nu}}) &\mapsto& (M,Y^{\hat{\nu}}_K),\\
\end{eqnarray*}
with the obvious definitions on morphisms.

The following theorem is a special case of the main results proved in
\cite{BDM}:

\begin{thm} [\cite{BDM}] \label{bdm-theorem}
The functors $F_{\hat{\nu}}$ and $G_{\hat{\nu}}$ have the properties
mentioned above.  Furthermore, $F_{\hat{\nu}} \circ G_{\hat{\nu}} =
id_{\mathcal{ C}^{\hat{\nu}}(V_L)}$ and $G_{\hat{\nu}} \circ
F_{{\hat{\nu}}} = id_{\mathcal{ C}(V_K)}$.  In particular, the
categories $\mathcal{ C}(V_K)$ and $\mathcal{ C}^{\hat{\nu}}(V_L)$ are
isomorphic, as are the subcategories of irreducible objects.
\end{thm}

\section{Realizing the ``space-time" construction of a $\hat{\nu}$-twisted
$V_L$-module as a $V_K$-module}\label{compare-section}

Let $(V_L^T, Y^{\hat{\nu}})$ be the $\hat{\nu}$-twisted $V_L$-module
constructed in Section \ref{spacetime-section} following \cite{L1}, \cite{FLM2}
and
\cite{DL}.  Then by Theorem \ref{bdm-theorem}, $(V_L^T, Y^{\hat{\nu}})$ is
isomorphic to some $\hat{\nu}$-twisted $V_L$-module $(M,Y_{\hat{\nu}})$
constructed via the method of Section \ref{bdm-section}, so that
$M$ is a $V_K$-module and $Y_{\hat{\nu}}$ is the twisted vertex
operator map defined by (\ref{bdm-twisted1}) and (\ref{bdm-twisted2}).
That is, $G_{\hat{\nu}} (V_L^T)$ must be a
$V_K$-module, and since $V_L^T$ is irreducible as a
$\hat{\nu}$-twisted module, $G_{\hat{\nu}}(V_L^T)$ must be an irreducible
$V_K$-module.  We shall write the conformal element of the vertex
operator algebra $V_K$ as $\omega_K$ and the corresponding Virasoro
algebra operators simply as $L(n)$, and we shall keep the same notation
$Y$ and ${\bf 1}$ for the vertex operator map and the vacuum vector of
$V_K$; that is, $V_K = (V_K, Y, {\bf 1}, \omega_K)$.  The central
charge of $V_K$ is $d$.

Irreducible modules for lattice vertex operator algebras are
classified as follows (\cite{D1}, \cite{DLiM1}; cf. \cite{LL}):
Let ${\mathcal L}$ be a positive definite even lattice and let
${\mathcal L}^*$ be the dual lattice to ${\mathcal L}$.  The
irreducible $V_{\mathcal L}$-modules are parametrized up to
equivalence by ${\mathcal L}^* / {\mathcal L}$, and in fact, each is
isomorphic to a ``coset module'' $V_{\beta + {\mathcal L} }$ for some
$\beta \in {\mathcal L}^*$.

Thus $G_{\hat{\nu}} (V_L^T)$ is isomorphic to $V_{\beta + K}$ for some coset
$\beta + K \in K^*/K$.  Writing $\h_K = K \otimes_\Z \C$ to
distinguish {}from $\h = L\otimes_\Z \C$, we have
that $V_{\beta +K} \simeq S(\hat{\h}^-_K) \otimes \C[\beta + K]$
linearly.  The grading of $V_{\beta +K}$ is given by weights with
${\rm wt} \; \alpha(-n) = n$ for $\alpha \in K$, $n \in \Z_+$, and
${\rm wt} \; e^{\beta + \alpha} = \frac{1}{2} \langle \beta + \alpha,
\beta + \alpha \rangle$ for $\alpha \in K$.

\begin{thm}\label{untwisted-isomorphism}
As a $V_K$-module, $G_{\hat{\nu}} (V_L^T)$ is isomorphic to $V_K$.
\end{thm}

\begin{proof}
As a $V_K$-module, $G_{\hat{\nu}}(V_L^T)$ is given by the space $V_L^T =
S[\nu]\otimes U_T \simeq S(\hat{\h}[\nu]^-) \otimes \C[P_0L] $ and the
vertex operators
\[ Y_K^{\hat{\nu}}(u,x) = Y^{\hat{\nu}}((\mathcal{E}_f(x)^{-1} u)^1, x^k) \]
for $u \in V_K$.

To determine the graded dimension of the $V_K$-module $G_{\hat{\nu}}(V_L^T)$,
we first observe that
\[ Y_K^{\hat{\nu}}(\omega_K,x) = \sum_{n \in \Z} L_K^{\hat{\nu}}(n) x^{-n-2} =
Y^{\hat{\nu}}((\mathcal{E}_f(x)^{-1} \omega_K)^1, x^k) . \]
We calculate $\mathcal{E}_f (x)^{-1} \omega_K$ by noticing that since
$L(j) {\bf 1} = 0$ for $j \geq -1$ and $\omega_K = L(-2) {\bf 1}$,
we have that if
$j \geq 1$, then $L(j) \omega_K = \frac{j^3 - j}{12} \delta_{j-2,0}
c{\bf 1}$.  Thus recalling that $a_2 = (k^2 - 1)/12$, we
have
\begin{eqnarray*}
\mathcal{E}_f (x)^{-1} \omega_K &=& x^{(k-1) L(0)} k^{L(0)} \exp \Biggl(
-\sum_{j \in \Z_+} a_j x^{- j} L(j) \Biggr) \cdot \omega_K \\
&=& x^{(k-1)L(0)} k^{L(0)} \left(\omega_K - a_2 x^{-2} \frac{1}{2} d {\bf
1}\right) \\
&=& x^{2k-2} k^2  \omega_K - \frac{(k^2-1)d}{24x^2} {\bf 1}.
\end{eqnarray*}
Thus
\begin{eqnarray*}
Y_K^{\hat{\nu}}(\omega_K,x) &=& Y^{\hat{\nu}} \Bigl( \bigl(x^{2k-2} k^2 \omega_K
- \frac{(k^2-1)
d}{24x^2} {\bf 1} \bigr)^1, x^k \Bigr) \\
&=& x^{2k-2} k^2 Y^{\hat{\nu}}(\omega_K^1, x^k) - \frac{(k^2-1)d}{24x^2} .\\
\end{eqnarray*}

Next we calculate $e^{\Delta_x} \cdot \omega_K^1$ by recalling that
$\omega_K = \frac{1}{2} \sum_{j=1}^d h_j (-1) h_j(-1) {\bf 1}$ where
$\{h_1, \dots, h_d \}$ is an orthonormal basis for $\h_K = K
\otimes_{\Z} \C$.   Note that since $\{h_1,\dots,h_d\}$ is an orthonormal basis
for $\h_K$, then
\begin{eqnarray}\label{ortho-basis}
\{\beta_1, \dots, \beta_{kd} \} &=& \{(h_1, 0,0,\dots,0), \dots,
(h_d,0,0,\dots,0), (0, h_1,0,\dots,0),
\nonumber \\
& & \dots, (0,h_d,0,\dots,0), \dots, (0,0,\dots,0,h_1),\dots, (0,0,\dots,0,h_d)
\} \nonumber \\
&=& \{h_j^p \ | \ j = 1,\dots,d \ \mathrm{and} \  p = 1,\dots, k \}
\end{eqnarray}
is an orthonormal basis for $\h = L \otimes_\Z \C$.  Thus using (\ref{c_110}),
we have
\begin{eqnarray*}
\lefteqn{ \Delta_x \cdot \omega_K^1 }\\
&=& \sum_{m,n\ge 0} \sum_{r=0}^{k-1} \sum^{kd}_{j=1} c_{mnr} (\nu^{-r}
\beta_j)(m) \beta_j (n) x^{-m-n} \cdot \frac{1}{2} \sum_{s=1}^d (h_s (-1)
h_s(-1) {\bf 1})^1 \\
&=& \sum_{m,n\ge 0} \sum_{r=0}^{k-1} \sum^{d}_{j=1} \sum_{p = 1}^k c_{mnr}
(\hat{\nu}^{-r} h^p_j)(m) h_j^p (n) x^{-m-n} \cdot \frac{1}{2} \sum_{s=1}^d (h_s
(-1) h_s(-1) {\bf 1})^1 \\
&=& \sum_{m,n\ge 0} \sum_{r=0}^{k-1} \sum^{d}_{j=1} \sum_{p=1}^k c_{mnr}
h_j^{r+p}
(m) h_j^p (n)  x^{-m-n} \cdot \frac{1}{2} \sum_{s=1}^d (h_s (-1) h_s(-1) {\bf
1})^1  \\
\end{eqnarray*}
\begin{eqnarray*}
&=& \frac{1}{2} \sum^{d}_{j=1} c_{110} \Bigl(h_j^{1} (1) h_j^1 (1) \Bigr)
x^{-2} \cdot (h_j (-1) h_j(-1) {\bf 1})^1 \hspace{1.5in}\\
&=& \sum^{d}_{j=1}  \frac{k^2-1}{24k^2}  \langle h_j , h_j \rangle^2 {\bf
1}\otimes {\bf
1} \otimes \cdots \otimes {\bf 1} x^{-2}\\
&=& \sum^{d}_{j=1}  \frac{k^2-1}{24k^2}  {\bf 1}\otimes {\bf 1} \otimes \cdots
\otimes{\bf 1} x^{-2}\\
&=&  \frac{(k^2-1)d}{24k^2}  {\bf 1}\otimes {\bf 1} \otimes \cdots \otimes{\bf
1} x^{-2} .
\end{eqnarray*}
and
\begin{equation}
e^{\Delta_x} \omega_K^1 = \omega_K^1 + \frac{(k^2-1)d}{24k^2} {\bf 1} \otimes
{\bf 1} \otimes \cdots \otimes {\bf 1} x^{-2} .
\end{equation}

Therefore,
\begin{eqnarray*}
\lefteqn{Y_K^{\hat{\nu}}(\omega_K,x) }\\
&=& x^{2k-2} k^2 Y^{\hat{\nu}}(\omega_K^1, x^k) - \frac{(k^2-1)d}{24x^2} \\
&=& x^{2k-2} k^2 W(e^{\Delta_{x^k}} \omega_K^1, x^k) - \frac{(k^2-1)d}{24x^2} \\
&=& x^{2k-2} k^2 W(\omega_K^1, x^k) + x^{2k-2} k^2  \frac{(k^2-1)d}{24k^2}
x^{-2k} - \frac{(k^2-1)d}{24x^2} \\
&=& x^{2k-2} \frac{k^2}{2} W\left(\sum_{j=1}^d ( h_j(-1) h_j(-1){\bf 1} )^1,
x^k\right) \\
&=& x^{2k-2} \frac{k^2}{2} \sum_{j=1}^d \ _\circ^\circ (h_j^1)^T(x^k)
(h_j^1)^T(x^k) _\circ^\circ\\
&=& x^{2k-2} \frac{k^2}{2} \sum_{j=1}^d \sum_{m,n \in \frac{1}{k} \Z} \
_\circ^\circ (h_j^1)^T(m) (h_j^1)^T(n) _\circ^\circ \, x^{-km-kn-2k} \\
&=& \frac{k^2}{2} \sum_{j=1}^d \sum_{m,n \in \frac{1}{k} \Z} \ _\circ^\circ
(h_j^1)^T(m) (h_j^1)^T(n)   _\circ^\circ  \, x^{-km-kn-2}.
\end{eqnarray*}
Thus
\begin{eqnarray*}
L_K^{\hat{\nu}}(0) &=& {\rm Res}_x x \biggl( \frac{k^2}{2} \sum_{j=1}^d
\sum_{m,n \in
\frac{1}{k} \Z} \ _\circ^\circ (h_j^1)^T(m) (h_j^1)^T(n) _\circ^\circ \,
x^{-km-kn-2} \biggr)\\
&=& \frac{k^2}{2} \sum_{j=1}^d \sum_{n \in \frac{1}{k} \Z} \ _\circ^\circ
(h_j^1)^T(-n) (h_j^1)^T(n) _\circ^\circ \\
&=& \frac{k^2}{2} \sum_{j=1}^d \sum_{n \in \frac{1}{k} \Z}
(h_j^1)_{(-k|n|)}(-|n|) (h_j^1)_{(k|n|)}(|n|) .
\end{eqnarray*}

We want to compare $L_K^{\hat{\nu}}(0)$ to $L^{\hat{\nu}}(0)$ given by
(\ref{L-twisted}).   To
do this,
we note that  for $n,p = 0,\dots,k-1$ and $h \in \h$, we have
\begin{equation}
(\nu^p h)_{(n)} =  \eta^{np}  h_{(n)}.
\end{equation}
Thus recalling (\ref{L-twisted}), we have
\begin{eqnarray*}
L^{\hat{\nu}}(0) &=&  \frac{1}{2} \sum_{j = 1}^{kd} \sum_{n \in \frac{1}{k}
\mathbb{Z}}
\beta_j^T(-|n|) \beta_j^T(|n|) + \frac{(k^2 - 1)d}{24k}  \\
&=& \frac{1}{2} \sum_{j = 1}^{d} \sum_{n \in \frac{1}{k} \Z} \sum_{p=1}^{k}
(h_j^p)_{(-k|n|)}(-|n|) (h_j^p)_{(k|n|)} (|n|) + \frac{(k^2 - 1)d}{24k} \\
&=& \frac{1}{2} \sum_{j = 1}^{d} \sum_{n \in \frac{1}{k} \Z} \sum_{p=1}^{k}
(\hat{\nu}^{-p+1} h_j^1)_{(-k|n|)}(-|n|) (\hat{\nu}^{-p+1}h_j^1)_{(k|n|)} (|n|)
+ \frac{(k^2
- 1)d}{24k} \\
&=& \frac{1}{2} \sum_{j = 1}^{d} \sum_{n \in \frac{1}{k} \Z} \sum_{p=1}^{k}
\eta^{-k|n| (-p+1)} ( h_j^1)_{(-k|n|)}(-|n|) \eta^{k|n| (-p+1) }
(h_j^1)_{(k|n|)} (|n|) \\
& & \quad + \frac{(k^2 - 1)d}{24k} \\
&=& \frac{k}{2} \sum_{j = 1}^{d} \sum_{n \in \frac{1}{k} \Z} (
h_j^1)_{(-k|n|)}(-|n|) (h_j^1)_{(k|n|)} (|n|) + \frac{(k^2 - 1)d}{24k} \\
&=& \frac{1}{k} L_K^{\hat{\nu}}(0) + \frac{(k^2 - 1)d}{24k} .
\end{eqnarray*}
In other words,
\begin{equation}\label{comparing-L(0)'s}
L_K^{\hat{\nu}}(0) = kL^{\hat{\nu}}(0) -\frac{(k^2-1)d}{24},
\end{equation}
which gives the natural grading on the space $V_L^T$ viewed as the $V_K$-module
$G_{\hat{\nu}}(V_L^T)$.  Thus by (\ref{character-twisted}) we have
\begin{eqnarray}
{\rm dim}_* G_{\hat{\nu}}( V_L^T) &=& {\rm tr}_{G_{\hat{\nu}}(V_L^T)} \;
q^{L_K^{\hat{\nu}}(0) - d/24} \\
&=& {\rm tr}_{V_L^T} \; q^{kL^{\hat{\nu}}(0) - (k^2-1)d/24 - d/24} \nonumber \\
&=& {\rm tr}_{V_L^T} \; q^{k(L^{\hat{\nu}}(0) - kd/24)} \nonumber \\
&=& \left. {\rm dim}_* V_L^T \right|_{q = q^k} \nonumber \\
&=& q^{-d/24} \left(\sum_{\alpha \in K} q^{\langle \alpha , \alpha \rangle /
2} \right) \left( \prod_{n \in \Z_+} (1 - q^n)^{-d} \right) . \nonumber
\end{eqnarray}
But the graded dimension of $V_K$ is given by
\begin{eqnarray}\label{dimension-equality}
{\rm dim}_* V_K &=& \frac{\Theta_K(q)}{\eta(q)^d} \\
&=& \left(\sum_{\alpha \in K} q^{ \langle \alpha , \alpha \rangle /2} \right)
q^{-d/24} \prod_{n \in \Z_+} (1 - q^n)^{-d},  \nonumber
\end{eqnarray}
so that
\begin{equation}
{\rm dim}_* G_{\hat{\nu}}( V_L^T) = {\rm dim}_* V_K .
\end{equation}
On the other hand, the graded dimension of the coset module $V_{\beta +K}$ with
$\beta \notin K$ is given by
\begin{eqnarray*}
{\rm dim}_* V_{\beta + K} &=& \left(\sum_{\alpha \in K} q^{ \langle \alpha ,
\alpha \rangle /2 + \langle \beta, \beta \rangle/2 + \langle \alpha, \beta
\rangle} \right) q^{-d/24} \prod_{n \in \Z_+} (1 - q^n)^{-d},
\end{eqnarray*}
which is of the form
\begin{eqnarray*}
{\rm dim}_* V_{\beta + K} &=& q^{-d/24} ( q^m + \cdots)
\end{eqnarray*}
for some $m \in \mathbb{Z}$, $m \neq 0$. Since ${\rm dim}_* V_K = q^{-d/24}
(1+ \cdots)$ we have that ${\rm dim}_* V_K \neq {\rm dim}_* V_{\beta + K}$ for
$K
\neq \beta + K \in K^*/K$, so that $V_K$ is the unique irreducible
$V_K$-module with graded dimension given by (\ref{dimension-equality}).
Therefore, $G_{\hat{\nu}}(V_L^T)$ and $V_K$ are isomorphic as $V_K$-modules.
\end{proof}

By Theorems \ref{bdm-theorem} and \ref {untwisted-isomorphism} we have the
following main result:

\begin{thm}\label{twisted-isomorphism}
The $\hat{\nu}$-twisted $V_K^{\otimes k}$-modules
$(V_L^T , Y^{\hat{\nu}})$ and $(V_K, Y_{\hat{\nu}})$ are isomorphic.
\end{thm}

In other words, the $\hat{\nu}$-twisted $V_L$-module constructed via the
``space-time" construction of Section \ref{spacetime-section} is isomorphic to
the $\hat{\nu}$-twisted $V_L$-module obtained {}from the $V_K$-module
$V_K$ using the ``worldsheet" construction of Section \ref{bdm-section}.

\section{An explicit determination of the isomorphism between the twisted
modules arising {}from the ``space-time" and the ``worldsheet"
constructions}\label{explicit-section}

We explicitly construct an isomorphism (necessarily unique up to
nonzero scalar multiple) given by Theorem \ref{twisted-isomorphism}.
This illuminates the correspondence between the two very different
twisted vertex operator maps.

Theorem \ref{twisted-isomorphism} gives the existence of a linear isomorphism
\begin{equation}
\mathcal{F}: V_L^T \longrightarrow V_K
\end{equation}
satisfying
\begin{equation}\label{iso-condition}
Y_{\hat{\nu}}(u,x) \mathcal{F}(v) = \mathcal{F}(Y^{\hat{\nu}}(u,x)v)
\end{equation}
for $u \in V_L \simeq V_K^{\otimes k}$ and $v \in V_L^T$.  In the next
theorem, we give the construction of this isomorphism $\mathcal{F}$ of
$V_K^{\otimes k}$-twisted modules, normalized so that
$\mathcal{F}(\mathbf{1}) = \mathbf{1}$.

\begin{thm}\label{iso-construction}
The normalized isomorphism $\mathcal{F} : V_L^T \longrightarrow V_K $ is the
unique linear map {}from $V_L^T$ to $V_K$ such that
\begin{equation}\label{main-iso}
\mathcal{F} \circ (\alpha,0,\dots,0)^T (x)\circ \mathcal{F}^{-1}  = \frac{1}{k}
x^{1/k-1} \alpha
(x^{1/k})
\end{equation}
for $\alpha \in K$ and such that $\mathcal{F}$ on the group algebra component
of $V_L^T$ is the isomorphism $\C[P_0L] \simeq \C[K]$ given by extension of
the isomorphism of $P_0 L \simeq K$,  $\frac{1}{k} (\alpha,\alpha,\dots,\alpha)
\mapsto
\alpha$, as in Remark \ref{lattice-quotient-remark}.  Furthermore,  for
$(\alpha_1, \dots, \alpha_k) \in L$, we have
\begin{equation}
\mathcal{F} \circ (\alpha_1, \dots, \alpha_k)^T(x) \circ \mathcal{F}^{-1} =
\frac{1}{k} x^{1/k-1} \sum_{j=1}^{k} \eta^{j-1} \alpha_j (\eta^{j-1} x^{1/k}).
\end{equation}
\end{thm}

\begin{proof} Let $\alpha \in K$.  Then
\begin{eqnarray*}
\lefteqn{\Delta_x \cdot (\alpha(-1)\iota(1) )^1 } \\
&=& \Delta_x \cdot (\alpha(-1)\iota(1) \otimes {\bf 1} \otimes \cdots \otimes
{\bf 1} ) \\
&=& \sum_{m,n\ge 0} \sum_{r=0}^{k-1} \sum^{d}_{j=1} \sum_{p=1}^k c_{mnr}
h_j^{r+p} (m) h_j^p (n) x^{-m-n} \cdot (\alpha (-1)\iota(1))^1  \\
&=& 0,
\end{eqnarray*}
and thus
\begin{eqnarray*}
Y^{\hat{\nu}}( ( \alpha(-1)\iota(1))^1 , x)&=& W(e^{\Delta_x}
(\alpha(-1)\iota(1))^1, x) \\
&=& W( (\alpha(-1)\iota(1) )^1, x) \\
&=& (\alpha,0,\dots, 0)^T (x)  .
\end{eqnarray*}
We also have
\begin{eqnarray*}
\mathcal{E}_f(x^{1/k}) \alpha(-1)\iota(1) &=& \exp \Biggl( \sum_{j \in \Z_+}
a_j x^{- j/k} L(j)  \Biggr) k^{-L(0)} x^{\left( 1/k -1 \right) L(0)} \alpha(-1)
\iota (1) \\
&=& \frac{1}{k} x^{1/k - 1} \alpha(-1) \iota(1),
\end{eqnarray*}
and thus
\begin{eqnarray*}
Y_{\hat{\nu}}( (\alpha(-1)\iota(1))^1, x)  &=& Y(\mathcal{E}_f(x^{1/k})
\alpha(-1)\iota(1), x^{1/k}) \\
&=& \frac{1}{k} x^{1/k-1} Y( \alpha(-1)\iota(1), x^{1/k}) \\
&=& \frac{1}{k} x^{1/k-1} \alpha(x^{1/k})  .
\end{eqnarray*}
Suppose that $\mathcal{F}$ is an isomorphism of  $\hat{\nu}$-twisted
$V_K^{\otimes
k}$-modules {}from $V_L^T$ to $V_K$; by Theorem \ref{twisted-isomorphism},
$\mathcal{F}$ exists.  Then since $\mathcal{F}$
must satisfy (\ref{iso-condition}), we have
\begin{eqnarray}\label{first-part-proof}
\mathcal{F} \circ (\alpha,0,\dots, 0)^T(x) \circ \mathcal{F}^{-1} &=&
\mathcal{F} \circ Y^{\hat{\nu}}((\alpha(-1)\iota(1))^1 ,x) \circ
\mathcal{F}^{-1}  \\
&=& Y_{\hat{\nu}} (\alpha(-1)\iota(1))^1  ,x)  \nonumber \\
&=& \frac{1}{k} x^{1/k-1} \alpha (x^{1/k}) \nonumber,
\end{eqnarray}
proving (\ref{main-iso}).

Let
\begin{eqnarray}
e : L = K\oplus K \oplus \cdots \oplus K & \longrightarrow & \hat{L}\\
(\alpha_1, \dots, \alpha_k ) & \mapsto & e_{(\alpha_1, \dots, \alpha_k)}
\nonumber
\end{eqnarray}
be a section of $\hat{L}$.  This choice of section allows us to identify
$\mathbb{C}\{L\}$ with the group algebra $\mathbb{C}[L]$ by the linear
isomorphism
\begin{eqnarray}
\mathbb{C}[L] & \longrightarrow & \mathbb{C}\{L\} \\
e^{(\alpha_1, \dots, \alpha_k)} & \mapsto & \iota(e_{(\alpha_1, \dots, \alpha_k
)}) \nonumber
\end{eqnarray}
for $\alpha_1, \dots, \alpha_k \in K$.  Without confusion, we use
the same notation for a section of $\hat{K}$. Then using the
identification of $V_L$ and $V_K^{\otimes k}$, we have
\begin{eqnarray*}
\lefteqn{Y^{\hat{\nu}}((e_\alpha)^1 ,x){\bf 1} }\\
&=& Y^{\hat{\nu}}( e_{(\alpha,0, \dots, 0)} ,x){\bf 1} \\
&=& k^{-\langle (\alpha,0,\dots,0) ,(\alpha, 0,\dots,0) \rangle /2}
\sigma((\alpha, 0, \dots, 0) \  _\circ^\circ \,
e^{\int( (\alpha, 0,\dots,0)^T (x)- (\alpha, 0,\dots,0)^T(0)x^{-1})} \cdot \\
& & \quad \cdot e_{(\alpha, 0,\dots,0)} x^{(\alpha, 0,\dots,0)_{(0)}+\langle
(\alpha, 0,\dots,0)_{(0)} ,(\alpha, 0,\dots,0) _{(0)}\rangle /2-\langle
(\alpha, 0,\dots,0),(\alpha, 0,\dots,0) \rangle /2} \ _\circ^\circ \, {\bf 1}
\\
&=& k^{-\langle \alpha , \alpha \rangle /2} \exp \biggl( \sum_{n \in
\frac{1}{k} \Z_+} \frac{ (\alpha,0,\dots,0)_{(-kn)} (-n)}{n} x^n  \biggr) \cdot
\\
& & \quad \cdot \exp\biggl( \sum_{n \in \frac{1}{k} \Z_+} \frac{
(\alpha,0,\dots,0)_{(kn)} (n)}{-n} x^{-n}\biggr)  e_{(\alpha, 0,\dots,0)} \cdot
\\
& & \quad \cdot x^{(\alpha,\alpha, \dots, \alpha)/k+\langle
(\alpha,\alpha,\dots,\alpha) ,(\alpha,\alpha, \dots, \alpha) \rangle
/(2k^2)-\langle \alpha,\alpha \rangle /2} {\bf 1} \\
&=& k^{-\langle \alpha,\alpha \rangle /2} \exp \biggl( \sum_{n \in \frac{1}{k}
\Z_+} \frac{ (\alpha,0,\dots, 0)_{(-kn)} (-n)}{n} x^n \biggr) x^{\langle
\alpha,\alpha \rangle/(2k) -\langle \alpha,\alpha \rangle /2} \iota(
e_{(\alpha, 0, \dots, 0)}) \\
&=& k^{-\langle \alpha,\alpha \rangle /2} \exp \biggl( \sum_{n \in \frac{1}{k}
\Z_+} \frac{ (\alpha,0,\dots, 0)_{(-kn)} (-n)}{n} x^n \biggr) x^{(1-k)\langle
\alpha,\alpha \rangle/(2k) } \iota(e_{(\alpha, 0,\dots, 0)})_{(0)},
\end{eqnarray*}
whereas
\begin{eqnarray*}
\lefteqn{Y_{\hat{\nu}}( (e_\alpha)^1 ,x) {\bf 1}} \\
&=& Y(\mathcal{E}_f (x^{1/k}) \iota(e_\alpha) ,x^{1/k}) {\bf 1} \\
&=& Y\Biggl( \exp \biggl( \sum_{j \in \Z_+} a_j x^{- j/k} L(j) \biggr)
k^{-L(0)} x^{\left( 1/k -1 \right) L(0)} \iota(e_\alpha) ,x^{1/k} \Biggr) {\bf
1} \hspace{.5in} \\
&=& Y( k^{-\langle \alpha,\alpha \rangle/2 } x^{(1/k-1) \langle \alpha, \alpha
\rangle/2 } e_\alpha ,x^{1/k}) {\bf 1} \\
&=& k^{-\langle \alpha,\alpha \rangle/2 } x^{ (1-k) \langle \alpha, \alpha
\rangle/(2k) } \ _\circ^\circ \, e^{\int(\alpha(x^{1/k})- \alpha(0)x^{-1/k})}
e_\alpha x^{\alpha /k} \ _\circ^\circ \, {\bf 1} \\
&=& k^{-\langle \alpha,\alpha \rangle/2 } x^{ (1-k)\langle \alpha, \alpha
\rangle/(2k) } \exp \Bigl( \sum_{n \in \Z_+} \frac{\alpha(-n)}{n} x^{n/k}
\Bigr) \exp \Bigl( \sum_{n \in \Z_+} \frac{\alpha(n)}{-n} x^{-n/k} \Bigr)
 x^{\alpha/k} e_\alpha {\bf 1} \\
&=& k^{-\langle \alpha,\alpha \rangle/2 } x^{ (1-k)\langle \alpha, \alpha
\rangle/(2k) } \exp \Bigl( \sum_{n \in \Z_+} \frac{\alpha(-n)}{n} x^{n/k}
\Bigr) \iota(e_\alpha) .
\end{eqnarray*}
By (\ref{first-part-proof}), we have
\begin{equation*}
\mathcal{F} \circ (\alpha,0,\dots,0)_{(-kn)} (-n)  = \frac{1}{k} \alpha(-kn)
\circ \mathcal{F},
\end{equation*}
for $n \in \frac{1}{k} \Z_+$, and thus
\begin{eqnarray*}
\lefteqn{k^{-\langle \alpha,\alpha \rangle /2} x^{(1-k)\langle \alpha,\alpha
\rangle/(2k)} \exp \biggl( \sum_{n \in \Z_+} \frac{ \alpha (-n)}{n} x^{n /k}
\biggr) \mathcal{F} (\iota(e_{(\alpha, 0,\dots,0)})_{(0)} ) }\\
&=& \mathcal{F} \Biggl( k^{-\langle \alpha,\alpha \rangle /2}x^{(1-k)\langle
\alpha,\alpha \rangle/(2k)} \exp \biggl( \sum_{n \in \Z_+} \frac{ k (\alpha,
0,\dots,0)_{(-n)} (-n/k)}{n} x^n \biggr) \cdot \\
& & \hspace{3.4in} \cdot \iota( e_{(\alpha, 0,\dots,0)})_{(0)} \Biggr) \\
&=& \mathcal{F} \Biggl( k^{-\langle \alpha,\alpha \rangle /2}x^{(1-k)\langle
\alpha,\alpha \rangle/(2k)} \exp \biggl( \sum_{n \in \frac{1}{k} \Z_+} \frac{
(\alpha, 0,\dots,0)_{(-kn)} (-n)}{n} x^n \biggr) \cdot \\
& & \hspace{3.4in} \cdot \iota( e_{(\alpha, 0,\dots,0)})_{(0)} \Biggr) \\
&=& \mathcal{F}( Y^{\hat{\nu}}((e_\alpha)^1 ,x){\bf 1}) \\
&=& Y_{\hat{\nu}}((e_\alpha)^1 ,x) \mathcal{F}({\bf 1}) \\
&=& Y_{\hat{\nu}}((e_\alpha)^1 ,x) {\bf 1}\\
&=& k^{-\langle \alpha,\alpha \rangle/2 } x^{ (1-k)\langle \alpha, \alpha
\rangle/(2k) } \exp \Bigl( \sum_{n \in \Z_+} \frac{\alpha(-n)}{n} x^{n/k}
\Bigr) \iota(e_\alpha) ,
\end{eqnarray*}
implying that
\begin{equation}\label{second-part-proof}
\mathcal{F}(\iota(e_{(\alpha, 0,\dots,0)})_{(0)}) =
\mathcal{F}(\iota(e_{\frac{1}{k} (\alpha, \alpha,\dots,\alpha)})) =
\iota(e_\alpha).
\end{equation}

The isomorphism $\mathcal{F}$ is uniquely defined by (\ref{first-part-proof})
and (\ref{second-part-proof}).

For $j= 0,\dots,k-1$, {}from (\ref{first-part-proof}) we have
\begin{eqnarray*}
\lefteqn{\mathcal{F} \circ (\nu^{-j} (\alpha,0,\dots,0))^T (x) \circ
\mathcal{F}^{-1} }\\
&=& \sum_{n \in \frac{1}{k} \Z} \mathcal{F} \circ (\nu^{-j} (
\alpha, 0, \dots, 0))_{(kn)} (n) \circ \mathcal{F}^{-1} x^{-n-1} \\
&=& \sum_{n \in \frac{1}{k} \Z} \eta^{-knj}  \mathcal{F} \circ (
\alpha, 0, \dots, 0)_{(kn)} (n) \circ \mathcal{F}^{-1} x^{-n-1} \\
&=& \sum_{n \in \frac{1}{k} \Z} \frac{\eta^{-knj}}{k} \alpha (kn) x^{-n-1}
\\
&=& \sum_{n \in \Z} \frac{\eta^{-nj}}{k} \alpha (n) x^{-n/k-1} \\
&=& \frac{\eta^{j}}{k} x^{1/k-1} \sum_{n \in \Z} \alpha (n) (\eta^{j}
x^{1/k})^{-n-1} \\
&=& \frac{\eta^{j}}{k} x^{1/k-1} \alpha (\eta^{j} x^{1/k}) .
\end{eqnarray*}
Thus in general, for $(\alpha_1,\dots, \alpha_k) \in L$, we have
\begin{eqnarray*}
\mathcal{F} \circ (\alpha_1,\dots, \alpha_k)^T(x) \circ
\mathcal{F}^{-1} &=& \sum_{j=1}^{k} \mathcal{F} \circ \nu^{-j+1}
(\alpha_j,0,\dots,0)^T(x) \circ
\mathcal{F}^{-1} \\
&=& \sum_{j=1}^{k} \frac{\eta^{j-1}}{k} x^{1/k-1} \alpha_j (\eta^{j-1} x^{1/k})
\\
&=& \frac{1}{k} x^{1/k-1} \sum_{j=1}^{k}  \eta^{j-1}  \alpha_j (\eta^{j-1}
x^{1/k}) .
\end{eqnarray*}
\end{proof}

We also note that in our construction of $\hat{\nu}$-twisted
$V_L$-modules we have restricted $\nu$ to be the particular
permutation automorphism of $L$ which cyclicly permutes the direct sum
components $K$ of $L$ by the $k$-cycle $(1 \; 2 \cdots k)$.
However, the setting of $\hat{g}$-twisted $V_L$-modules makes sense for
$g$ an arbitrary permutation on $k$ letters acting on $L$.  It is easy to
extend these results to an arbitrary $k$-cycle permutation, since
any $k$-cycle is equal to $\mu \nu \mu^{-1}$ for some permutation $\mu$ and
$\nu = (1 \; 2 \cdots k)$.  The category of $\hat{\nu}$-twisted $V_L$-modules
$\mathcal{C}^{\hat{\nu}}(V_L)$ is isomorphic to the category of $\hat{\mu}
\hat{\nu} \hat{\mu} ^{-1}$-twisted $V_L$-modules $\mathcal{C}^{\hat{\mu}
\hat{\nu} \hat{\mu}^{-1}}(V_L)$ with the isomorphism given by
\begin{eqnarray}
H_{\hat{\mu}} : \mathcal{C}^{\hat{\nu}}(V_L) & \longrightarrow &
\mathcal{C}^{\hat{\mu} \hat{\nu}
\hat{\mu}^{-1}}(V_L) \\
(M,Y_{\hat{\nu}}) & \mapsto & (M, Y_{\hat{\mu} \hat{\nu} \hat{\mu}^{-1}})
\nonumber
\end{eqnarray}
where
\begin{equation}
Y_{\hat{\mu} \hat{\nu} \hat{\mu}^{-1}} (v,x) = Y_{\hat{\nu}} ( \hat{\mu} v, x)
\end{equation}
for $v \in V_L$ (cf. \cite{DLiM}, \cite{BDM}).  Thus our isomorphism
between $\hat{\nu}$-twisted $V_L$-modules extends to $\hat{g}$-twisted
$V_L$-modules for $g$ an arbitrary $k$-cycle.

For an arbitrary permutation on $k$ letters, $g$, we note that $g$ can
be written as a product of disjoint cycles $g=g_1\cdots g_p$ where the
order of $g_i$ is $k_i$ such that $\sum_{i}k_i=k$.  (Note that we are
including 1-cycles.)  Following \cite{BDM}, we further note that there exists
a permutation $\mu$ on $k$ letters satisfying $g= \mu g_1'\cdots g_p'
\mu^{-1}$ such that $g_i'$ is a $k_i$-cycle which permutes the numbers
\begin{equation}
\Bigl(\sum_{j = 1}^{i-1} k_j\Bigr) + 1, \; \Bigl(\sum_{j = 1}^{i-1} k_j\Bigr) +
2, \; \dots, \sum_{j = 1}^i k_j.
\end{equation}
We have already determined how to construct the $\hat{g_i'}$-twisted
$V^{\otimes k_i}$-module $V_{L_i}^T$ (where $L_i$ is the orthogonal
direct sum of $K_i$ copies of $K$) using the method of \cite{L1},
\cite{FLM2}, \cite{DL} and the $\hat{g_i'}$-twisted $V^{\otimes
k_i}$-module $V_K$ using the method of \cite{BDM}.  We have also
determined the isomorphisms $\mathcal{F}_i$ between these two
constructions $V_{L_i}^T$ and $V_K$.  {}From \cite{BDM}, we then have
the construction of the $\hat{g}$-twisted $V_L$-module $V_K^{\otimes
p}$, and putting the isomorphisms $\mathcal{F}_i$ together with the
isomorphism related to the conjugation $\mu$, we have an isomorphism
between the $\hat{g}$-twisted $V_L$-module $V_L^T$ of \cite{L1},
\cite{FLM2}, \cite{DL}, and the $\hat{g}$-twisted $V^{\otimes
k}$-module $V_K^{\otimes p}$ of \cite{BDM}.

Note that in the discussion above, we have not specified a unique
decomposition $g = \mu g_1'\cdots g_p' \mu^{-1}$ but rather have shown
how to construct the isomorphism of $\hat{g}$-twisted
$V^{\otimes k}$-modules for a given such (non-unique) decomposition
$g = \mu g_1'\cdots g_p' \mu^{-1}$.  However, for any decomposition $g = \mu
g_1'\cdots g_p' \mu^{-1}$, the resulting $\hat{g}$-twisted $V^{\otimes
k}$-module is isomorphic to that obtained {}from any other
decomposition.

\begin{rema}{\em
Finally, we comment that our isomorphisms between twisted modules also
carry over to the still more general case of lattice-cosets.  In
Section 10 (``Shifted vertex operators and their commutators'') of
\cite{L1}, the construction of the spaces $V_L^T$ recalled in Section
\ref{spacetime-section} above was in fact generalized to the following
setting: an arbitrary positive-definite even lattice $L$, an arbitrary
isometry $\nu$, and an arbitrary {\it coset} $L+\gamma$ of $L$, where
$\gamma$ is any element of the $\nu$-fixed vector subspace of the
ambient vector space spanned by $L$; this was a matter of ``shifting''
the original even-lattice construction in \cite{L1}.  For the
situation in which $\gamma$ is taken to be an element of the rational
span of $L$, this ``lattice-coset'' twisted-module construction,
including the enhancement of the structure given in \cite{FLM2},
\cite{L2} and \cite{DL}, is a special case of the theory carried out
in \cite{DL}.  In our setting in this paper, we may take $\gamma$ to
be any element of the dual lattice of the fixed-sublattice $K$ of $L$,
and we find that our isomorphism given in Theorem
\ref{twisted-isomorphism} generalizes to {\it all} the irreducible
$\hat{\nu}$-twisted $V_K^{\otimes k}$-modules and {\it all} the
irreducible $V_K$-modules.  In addition, this more general result of
course extends still further to arbitrary permutations, as described
above.}
\end{rema}

\end{document}